\documentclass[letterpaper,english]{smfart}
\usepackage{amsmath,amssymb}
\usepackage[pagebackref]{hyperref}

\renewcommand{\Re}{\operatorname{Re}}
\renewcommand{\Im}{\operatorname{Im}}

\newcommand{\cK}{{\mathcal K}}
\newcommand{\cM}{{\mathcal M}}
\newcommand{\cX}{{\mathcal X}}

\newcommand{\bbC}{{\mathbb C}}
\newcommand{\bbF}{{\mathbb F}}
\newcommand{\bbH}{{\mathbb H}}
\newcommand{\bbO}{{\mathbb O}}
\newcommand{\bbR}{{\mathbb R}}
\newcommand{\bbS}{{\mathbb S}}
\newcommand{\bbZ}{{\mathbb Z}}

\newcommand{\euso}{\operatorname{\mathfrak{so}}}
\newcommand{\euspin}{\operatorname{\mathfrak{spin}}}

\newcommand{\eusl}{\operatorname{\mathfrak{sl}}}

\newcommand{\eugl}{\operatorname{\mathfrak{gl}}}
\newcommand{\eug}{\operatorname{\mathfrak g}}
\newcommand{\euh}{\operatorname{\mathfrak h}}

\newcommand{\eup}{\operatorname{\mathfrak p}}
\newcommand{\eur}{\operatorname{\mathfrak r}}
\newcommand{\euk}{\operatorname{\mathfrak k}}

\newcommand{\SO}{\operatorname{SO}}\newcommand{\SL}{\operatorname{SL}}
\newcommand{\SU}{\operatorname{SU}}\newcommand{\Un}{\operatorname{U}}
\newcommand{\GL}{\operatorname{GL}}\newcommand{\Symp}{\operatorname{Sp}}
\newcommand{\tr}{\operatorname{tr}}\newcommand{\Or}{\operatorname{O}}
\newcommand{\Hom}{\operatorname{Hom}}
\newcommand{\I}{\operatorname{I}}
\newcommand{\G}{\operatorname{G}}\newcommand{\Pf}{\operatorname{Pf}}
\newcommand{\Spin}{\operatorname{Spin}}\newcommand{\Pin}{\operatorname{Pin}}
\newcommand{\Cl}{\operatorname{C\ell}}\newcommand{\Aut}{\operatorname{Aut}}
\newcommand{\Ric}{\operatorname{Ric}}\newcommand{\End}{\operatorname{End}}

\newcommand{\oneb}{\mathbf{1}}\newcommand{\zerob}{\boldsymbol{0}}
\newcommand{\bR}{\mathbf{R}}\newcommand{\bF}{\mathbf{F}}
\newcommand{\bH}{\mathbf{H}}
\newcommand{\sbold}{\mathbf{s}}
\newcommand{\hb}{\mathbf{h}}
\newcommand{\ab}{\mathbf{a}}
\newcommand{\qb}{\mathbf{q}}
\newcommand{\ub}{\mathbf{u}}
\newcommand{\xb}{\mathbf{x}}\newcommand{\yb}{\mathbf{y}}
\newcommand{\zb}{\mathbf{z}}

\newcommand{\bfomega}{\boldsymbol{\omega}}
\newcommand{\bftheta}{\boldsymbol{\theta}}
\newcommand{\bfeta}{\boldsymbol{\eta}}
\newcommand{\bfphi}{\boldsymbol{\phi}}

\newcommand{\bfTheta}{\boldsymbol{\Theta}}
\newcommand{\bfPhi}{\boldsymbol{\Phi}}
\newcommand{\bfsigma}{\boldsymbol{\sigma}}
\newcommand{\bfalpha}{\boldsymbol{\alpha}}

\newcommand{\p}{\partial}
\newcommand{\la}{\langle}\newcommand{\ra}{\rangle}

\newcommand{\ts}{\textstyle }

\newcommand{\w}{{\mathchoice{\,{\scriptstyle\wedge}\,}{{\scriptstyle\wedge}}
      {{\scriptscriptstyle\wedge}}{{\scriptscriptstyle\wedge}}}}

\author[R. Bryant]{Robert L. Bryant}
\address{Duke University Mathematics Department\\
         P.O. Box 90320\\
         Durham, NC 27708-0320}
\email{bryant@math.duke.edu}
\urladdr{http://www.math.duke.edu/\!\lower3pt\hbox{\symbol{'176}}bryant}

\title[Metrics with parallel spinors]
      {Pseudo-Riemannian metrics\\ 
       with parallel spinor fields\\
       and vanishing Ricci tensor }

\date{April 11, 2000}

\begin{document}

\begin{abstract}
   I will discuss geometry and normal forms for pseudo-Riemannian
   metrics with parallel spinor fields in some interesting dimensions.
   I also discuss the interaction of these conditions for parallel
   spinor fields with the Einstein equations.
\end{abstract}



\thanks{The research for this article was made possible by support 
        from the National Science Foundation through grant DMS-9870164
        and from Duke University.}

\maketitle

\setcounter{tocdepth}{2}
\tableofcontents

\section{Introduction}\label{sec:intro}

\subsection{Riemannian holonomy and parallel spinors}\label{ssec:holparspin}
The possible restricted holonomy groups of irreducible Riemannian 
manifolds have been known for some time now~\cite{mBerger55,rBr87,rBr96}.
The list of holonomy-irreducible types in dimension~$n$ 
that have nonzero parallel spinor fields is quite short:
The holonomy~$H$ of such a metric must be one of
\begin{itemize}
\item $H=\SU(m)$ (i.e., special K\"ahler metrics in dimension~$n=2m$);
\item $H=\Symp(m)$ (i.e., hyper-K\"ahler metrics in dimensions~$n=4m$);
\item $H=\G_2$ (when~$n=7$); or 
\item $H=\Spin(7)$ (when~$n=8$).
\end{itemize} 
In Cartan's sense, the local generality~\cite{rBr87,rBr96} 
of metrics with holonomy
\begin{itemize}
\item $H=\SU(m)$  ($n=2m$) is  $2$ functions of~$2m{-}1$ variables,
\item $H=\Symp(m)$  ($n=4m$) is $2m$ functions of~$2m{+}1$ variables,
\item $H=\G_2$  ($n=7$) is $6$ functions of~$6$ variables, and 
\item $H=\Spin(7)$  ($n=8$) is $12$ functions of~$7$ variables.
\end{itemize}
In each case, a metric with holonomy~$H$ has vanishing Ricci tensor. 

\subsection{Relations with physics}\label{ssec:relphy}
The existence of parallel spinor fields seems to account for much of 
the interest in metrics with special holonomy in mathematical physics,
since such spinor fields play a central role in supersymmetry.
In the case of string theory, $\SU(3)$, and lately, 
with the advent of $\cM$-theory, $\G_2$ (and possibly even $\Spin(7)$)
seem to be of interest.  I don't know much about these physical theories, 
so I will not attempt to discuss them.  

\subsection{Pseudo-Riemannian generalizations}\label{ssec:psdoRgen}
In the past few years, I have been asked by a number of physicists 
about the generality of pseudo-Riemannian metrics satisfying conditions 
having to do with parallel spinors and with solutions of the Einstein 
equations.  (In contrast to the Riemannian case,  an indecomposable 
pseudo-Riemannian metric can possess a parallel spinor field without 
being Einstein.)

For example, there seems to be some current interest in 
Lorentzian manifolds of type~$(10,1)$ having parallel spinor fields and
perhaps also having vanishing Ricci curvature, about which 
I will have more to say later in the article.

Recall~\cite{hW67,aBe87} that in the pseudo-Riemannian case, there is a 
distinction to be made between a metric being holonomy-irreducible 
(no parallel subbundles of the tangent bundle), being holonomy-indecomposable 
(no parallel splitting of the tangent bundle), and being indecomposable
(no local product decomposition of the metric).  (In the Riemannian case,
of course, these conditions are locally equivalent.)  The classification
of the holonomy-irreducible case proceeds much as in the 
positive definite case~\cite{rBr99}, but an indecomposable pseudo-Riemannian
metric need not be holonomy irreducible.  It is this difference
that makes classifying the possible pseudo-Riemannian metrics having
parallel spinor fields something of a challenge.    For a general
discussion of the differences, particularly the failure of the 
de~Rham splitting theorem, see~\cite{lBaI93,lBaI97}.  Also,
the results and examples in~\cite{aI96,aI99} are particularly
illuminating.

Now, quite a lot is known about the pseudo-Riemannian case when
the holonomy acts irreducibly.  For a general survey in this case, 
particularly regarding the existence of parallel spinor fields, 
see~\cite{hBiK99}.  Note that, in all of these cases, the Ricci
tensor vanishes.  This is not so when the holonomy acts reducibly. 
Already in dimension~$3$, Lorentzian metrics can have parallel spinor
fields without being Ricci-flat.  

An intriguing relationship between the condition for having a
parallel spinor and the Ricci equations came to my attention after
a discussion during a 1997 summer conference in Edinburgh with Ines Kath.
It had been known for a while~\cite{rBr87} that the metrics in
dimension~$7$ with holonomy~$\G_2$ depend locally on six functions of
six variables (modulo diffeomorphism).  Now, the condition of having
holonomy in~$\G_2$ is equivalent to the condition of having a parallel
spinor field.  I had also shown that the $(4,3)$-metrics with holonomy
$\G_2^*$ depend locally on six functions of six variables, and the
condition of having this holonomy in this group is the same as the
condition that the $(4,3)$-metric admit a non-null parallel spinor
field.  Ines Kath had noticed that the structure equations of a $(4,3)$
metric with a null parallel spinor field did not seem to imply that
the Ricci curvature vanished, and she wondered whether or not there
existed examples in which it did not.  After some analysis, 
I was able to show that there are indeed $(4,3)$-metrics with parallel 
spinor fields whose Ricci curvature is not zero and whose holonomy 
is equal to the full stabilizer of a null spinor.  
These metrics depend on three arbitrary
functions of seven variables.   However, a more intriguing result
is that, when one combines the condition of having a parallel
null spinor with the condition of being Ricci-flat, the $(4,3)$-metrics
with this property depend on six functions of six variables,
just as in the non-null case (where the vanishing of the Ricci
tensor is automatic).  

In any case, this and the questions from physicists motivates the general 
problem of determining the local generality of pseudo-Riemannian metrics 
with parallel spinors, with and without imposing the Ricci-flat condition.  
This article will attempt to describe some of what is known and give some 
new results, particularly in dimensions greater than~$6$.

Most of the normal forms that I describe for metrics
with parallel spinor fields of various different algebraic types 
are already known in the literature, or have been derived independently
by others.  (In particular, Kath~\cite{iKpurenote} has independently derived 
the normal forms for the split cases with a pure parallel spinor.)
What I find the most interesting is that, in every known case, 
the system of PDE given by the Ricci-flat condition is either in involution 
(in Cartan's sense) with the system of PDE that describe the $(p,q)$-metrics 
with a parallel spinor of given algebraic type or else follows as a
consequence (and so, in a manner of speaking, is trivially in 
involution with the parallel spinor field condition).  I have no
general proof that this is so in all cases, nor even a precise statement
as to how general the solutions should be, since this seems to depend
somewhat on the algebraic type of the parallel spinor.  What does
seem to be true in a large number of (though not all) cases, though, 
is that the local generality of the Ricci-flat $(p,q)$-metrics with
a parallel spinor of a given algebraic type seems to be largely independent
of the given algebraic type, echoing the situation for~$(4,3)$-metrics
mentioned above that first exhibited this phenomenon.

Since this article is mainly a discussion of cases, together with
an explicit working out of the standard moving frame methods 
and applications of Cartan-K\"ahler theory, I cannot claim a great 
deal of originality for the results.  Consequently, I do not state
the results in the form of theorems, lemmas, and propositions, but
instead discuss each case in turn.  The most significant
results are probably the descriptions of the generality of the 
Ricci-flat metrics with parallel spinors in the various cases.
Another possibly significant result is the description of the
$(10,1)$-metrics with a parallel null spinor field, since this
seems to be of interest in physics~\cite{jFF99}.

\section{Algebraic background on spinors}\label{sec:spinors}

All of the material in this section is classical.  I include it
to fix notation and for the sake of easy reference for the 
next section.  For more detail, the reader can consult~\cite{rHa90, LM89}.

\subsection{Notation}\label{ssec:notation}
The symbols~$\bbR$, $\bbC$, $\bbH$, and $\bbO$ denote, as usual,
the rings of real numbers, complex numbers, quaternions, and
octonions, respectively.  When~$\bbF$ is one of these rings, 
the notation~$\bbF(n)$ means the ring of $n$-by-$n$ matrices with
entries in~$\bbF$.  The notation~$\bbF^n$ will always denote
the space of \emph{column} vectors of height~$n$ with entries in~$\bbF$.
Vector spaces over~$\bbH$ will always be
regarded as having the scalar multiplication acting on the \emph{right}.
For an $m$-by-$n$ matrix~$a$ with entries in~$\bbC$ or~$\bbH$, the
notation~$a^*$ will denote its conjugate transpose.  When~$a$ has
entries in~$\bbR$, $a^*$ will simply denote the transpose of~$a$.

The notation~$\bbR^{p,q}$ denotes~$\bbR^{p+q}$ endowed with an
inner product of type~$(p,q)$.  The notation~$\bbC^{p,q}$ 
denotes~$\bbC^{p+q}$ endowed with an Hermitian inner product of 
type~$(p,q)$, with a similar interpretation of~$\bbH^{p,q}$, 
but the reader should keep in mind that a quaternion Hermitian inner 
product satisfies~$\la v  , w  q\ra = \la v  , w \ra q$
and $\la v  q , w \ra  = \bar q \la v  , w \ra $
for~$q\in\bbH$.

\subsection{Clifford algebras}\label{ssec:clifford}
The Clifford algebra~$\Cl(p,q)$ is the associative algebra 
generated by the elements of~$\bbR^{p,q}$ subject to the 
relations~$ v  w + w  v  = -2 v{\cdot}w \,\oneb$.  
This is a $\bbZ_2$-graded algebra, with
the even subalgebra~$\Cl^e(p,q)$ generated 
by the products~$vw$ for~$v,w\in\bbR^{p,q}$.

Because of the following formulae, valid for~$p,q\ge0$ 
(see~\cite{rHa90,LM89}),
\begin{equation}
\begin{split}
\Cl^e(p{+}1,q)&\simeq\Cl(p,q)\\
\Cl(p{+}1,q{+}1)&\simeq\Cl(p,q)\otimes\Cl(1,1)\\
\Cl(p{+}8,q)&\simeq \Cl(p,q)\otimes\Cl(8,0)\\
\Cl(p,q{+}1)&\simeq\Cl(q,p{+}1)\\
\end{split}
\end{equation}
all these algebras can be worked out from the table
\begin{equation}
\begin{split}
\Cl(0,1)&\simeq \bbR\oplus\bbR\\
\Cl(1,1)&\simeq \bbR(2)\\
\Cl(1,0)&\simeq \bbC\\
\Cl(2,0)&\simeq \bbH\\
\Cl(3,0)&\simeq \bbH\oplus\bbH\\
\Cl(4,0)&\simeq \bbH(2)\\
\Cl(5,0)&\simeq \bbC(4)\\
\Cl(6,0)&\simeq \bbR(8)\\
\Cl(7,0)&\simeq \bbR(8)\oplus\bbR(8)\\
\Cl(8,0)&\simeq \bbR(16).\\
\end{split}
\end{equation}
For example, $\Cl^e(p{+}1,p{+}1)\simeq\Cl(p,p{+}1)
\simeq\bbR(2^p)\oplus\bbR(2^p)$.

\subsection{$\Spin(p,q)$ and spinors}\label{ssec:spinpq}

By the defining relations, if~$v\cdot v\not=0$, then~$v\in\bbR^{p,q}$ is
a unit in~$\Cl(p,q)$ and, moreover, the twisted conjugation~$\rho(v):
\Cl(p,q)\to \Cl(p,q)$ defined on generators~$w\in\bbR^{p,q}$ 
by~$\rho(v)(w) = -vw v^{-1}$ preserves the generating subspace~$\bbR^{p,q}
\subset\Cl(p,q)$, acting as reflection in the hyperplane~$v^\perp
\subset\bbR^{p,q}$.  

The group~$\Pin(p,q)\subset\Cl(p,q)$ is the subgroup of the units 
in~$\Cl(p,q)$ generated by the elements~$v$ where~$v\cdot v = \pm1$ 
and the group~$\Spin(p,q)=\Pin(p,q)\cap\Cl^e(p,q)$ is 
the subgroup of the even Clifford algebra generated by
the products~$vw$, where~$ v {\cdot} v = w {\cdot} w =\pm1$.  

The map~$\rho$ defined above extends to a group homomorphism
$\rho:\Pin(p,q)\to\Or(p,q)$ that turns out to be a non-trivial double 
cover. The homomorphism~$\rho:\Spin(p,q)\to\SO(p,q)$ is also a 
non-trivial double cover.
 
The space of spinors~$\bbS^{p,q}$ is essentially
an irreducible~$\Cl(p,q)$-module, considered as
a representation of~$\Spin(p,q)$. 

When~$p{-}q\equiv 3\mod 4$, this definition is independent of 
which of the two possible irreducible~$\Cl(p,q)$ modules one 
uses in the construction.  

When~$p{-}q\equiv 0\mod 4$, the space~$\bbS^{p,q}$ is a reducible 
$\Spin(p,q)$-module, in fact, it can be written as a sum~$\bbS^{p,q}
=\bbS^{p,q}_+\oplus\bbS^{p,q}_-$ where~$\bbS^{p,q}_\pm$ are irreducible.
Action by an element of~$\Pin(p,q)$ not in~$\Spin(p,q)$ exchanges 
these two summands.

When~$p{-}q\equiv 1$ or~$2\mod 8$, the definition of~$\bbS^{p,q}$
as given above turns out to be the sum of two equivalent representations
of~$\Spin(p,q)$.  In this case, it is customary to redefine~$\bbS^{p,q}$
to be one of these two summands, so I do this without comment in the
rest of the article.

When~$q=0$, i.e., in the Euclidean case, I
will usually simplify the notation by writing~$\Cl(p)$, $\Spin(p)$,
and~$\bbS^p$ instead of~$\Cl(p,0)$, $\Spin(p,0)$, and~$\bbS^{p,0}$,
respectively.

\subsection{Orbits in the low dimensions}\label{ssec:orbitslow}

I will now describe the $\Spin(p,q)$-orbit structure of~$\bbS^{p,q}$ 
when~$p{+}q\le 6$.  This description made simpler by the fact that
there are several `exceptional isomorphisms' of Lie groups (as
discovered by Cartan) that reduce the problem to a series of classical
linear algebra problems.   

When~$p{+}q\le 1$, these groups are not particularly interesting and, 
since there is no holonomy in dimension~$1$ anyway, 
I will skip these cases.

\subsubsection{Dimension 2}\label{sssec:spin2}
Here there are two cases.

\paragraph{$\Spin(2)\simeq\Un(1)$}\label{par:spin20}
The action of~$\Spin(2)=\Un(1)$ 
on~$\bbS^2\simeq\bbC$ is the unit circle action
\begin{equation}
\lambda\cdot s = \lambda s\,.
\end{equation}
The orbits of~$\Spin(2)$ on~$\bbS^2 = \bbC$ are simply
the level sets of the squared norm, so all of the nonzero orbits
have the same stabilizer, namely, the identity.

Identifying~$\bbR^{2,0}$ with~$\bbC$, the action
of~$\Spin(2)$ on~$\bbR^{2,0}$ can be described as
\begin{equation}
\lambda\cdot v = \lambda^2\,v
\end{equation}
and the inner product is~$v\cdot v = |v|^2 = \bar v\,v$.

\paragraph{$\Spin(1,1)\simeq\bbR^*$}\label{par:spin11}
The action of~$\Spin(1,1)$ 
on~$\bbS^{1,1}\simeq\bbR\oplus\bbR$ is
\begin{equation}
\lambda\cdot(s_+,s_-) = (\lambda\,s_+, \lambda^{-1}\,s_-).
\end{equation}
There is an identification~$\bbR^{1,1}\simeq\bbR\oplus\bbR$
for which the action of~$\Spin(1,1)$ on~$\bbR^{1,1}$ has the description
\begin{equation}
\lambda\cdot(u,v) = (\lambda^2\,u, \lambda^{-2}\,v).
\end{equation}
and the inner product is~$(u,v)\cdot(u,v) = uv$. 
 
The nonzero orbits of~$\Spin(1,1)$ on~$\bbS^{1,1}$ are all of 
dimension~$1$ and have the same stabilizer, namely, the identity.

\subsubsection{Dimension 3}\label{sssec:spin3}
Again, there are two cases.

\paragraph{$\Spin(3)\simeq\Symp(1)$}\label{par:spin30} 
The action of~$\Spin(3)$ on~$\bbS^3\simeq\bbH$ is as 
quaternion multiplication:
\begin{equation}
A\cdot v = Av,
\end{equation}
where~$A$ and~$v$ are quaternions.
There are only two types of orbits, classified according to their
stabilizer types: Those of the point~$(0,0)$ and those of
the points~$(r,0)$, where~$r>0$ is a real number.  The stabilizer of
each nonzero element is trivial.

Identify~$\bbR^{3,0}$ with~$\Im\bbH$, so that the representation
of~$\Spin(3)$ on~$\bbR^{3,0}$ can be described as
\begin{equation}
A\cdot v = A\,v\,\overline{A}.
\end{equation}
and the inner product is~$v\cdot v = v\,\overline{v}$.

\paragraph{$\Spin(2,1)\simeq\SL(2,\bbR)$}\label{par:spin21}
The action of~$\Spin(2,1)$ on~$\bbS^{2,1}\simeq\bbR^2$ is as
the usual matrix multiplication:
\begin{equation}
A\cdot s  = A\,s.
\end{equation}
There are two $\Spin(2,1)$-orbits in~$\bbS^{2,1}$:
The orbit of the zero vector and then everything else.

Identify~$\bbR^{2,1}$ with the the space of
symmetric $2$-by-$2$ matrices, so that the representation
of~$\Spin(2,1)$ on~$\bbR^{2,1}$ can be described as
\begin{equation}
A\cdot v = A\,v\,A^*
\end{equation}
and the inner product is~$v\cdot v = -\det(v)$.  

There is an equivariant 
`spinor squaring' mapping~$\sigma:\bbS^{2,1}\to\bbR^{2,1}$ defined by
$\sigma(s) = s\,s^*$.   Its image is one nappe of the null cone
in~$\bbR^{2,1}$.

\subsubsection{Dimension 4}\label{sssec:spin4}
Now, there are three cases.

\paragraph{$\Spin(4)\simeq\Symp(1)\times\Symp(1)$}\label{par:spin40} 
The action of~$\Spin(4)$ on~$\bbS^{4}\simeq\bbH\oplus\bbH$ is
\begin{equation}
(A,B)\cdot(s_+,s_-) = (As_+, Bs_-).
\end{equation}
There are four types of spinor orbits (classified according to their
stabilizer types), those of the points~$(0,0)$, $(r_+,0)$,
$(0,r_-)$, and $(r_+,r_-)$, where~$r_\pm>0$ are real numbers.
Note that the stabilizer of a `generic' orbit (i.e., the fourth type)
is trivial.  Note that action by an element of~$\Pin(4)$ not in~$\Spin(4)$
exchanges the two summands and hence the two types of $3$-dimensional
orbits.

Under the identification~$\bbR^{4,0}\simeq\bbH$, the action
of~$\Spin(4)$ can be described as
\begin{equation}
(A,B)\cdot v = A\,v\,\overline{B}.
\end{equation}
and the inner product is~$v\cdot v = v\,\overline{v}$.

\paragraph{$\Spin(3,1)\simeq\SL(2,\bbC)$}\label{par:spin31}
The action of~$\Spin(3,1)$ on~$\bbS^{3,1}\simeq\bbC^2$ is
just
\begin{equation}
A \cdot s = As.
\end{equation}
In this case, there
are only two orbits, those of $0$ and $s$, where~$s\in\bbC^2$ is nonzero. 

Under the identification~$\bbR^{3,1}\simeq H_2(\bbC)$, the
Hermitian symmetric $2$-by-$2$ complex matrices, the action
of~$\Spin(3,1)$ can be described as
\begin{equation}
A \cdot v = A\,v\,A^*
\end{equation}
and the inner product is~$v\cdot v = -\det(v)$.  

There is an equivariant 
`spinor squaring' mapping~$\sigma:\bbS^{3,1}\to\bbR^{3,1}$ defined by
$\sigma(s) = s\,s^*$.   Its image is one nappe of the null cone
in~$\bbR^{3,1}$.  In relativity, this is referred to as the `forward
light cone'.

\paragraph{$\Spin(2,2)\simeq\SL(2,\bbR)\times\SL(2,\bbR)$}\label{par:spin22}
The action of~$\Spin(2,2)$ on~$\bbS^{2,2}\simeq\bbR^2\oplus\bbR^2$ is
\begin{equation}
(A,B)\cdot(s_+,s_-) = (As_+, Bs_-).
\end{equation}
There are four orbits of~$\Spin(2,2)$ on~$\bbS^{2,2}$, 
those of the points~$(0,0)$, $(s,0)$, $(0,s)$, and $(s,s)$, 
where~$s$ is any nonzero vector in~$\bbR^2$.  Note that action
by an element of~$\Pin(2,2)$ not in~$\Spin(2,2)$ exchanges the
two $2$-dimensional orbits.

Under the identification~$\bbR^{2,2}\simeq\bbR(2)$, the action
of~$\Spin(2,2)$ on~$\bbR^{2,2}$ can be described as
\begin{equation}
(A,B)\cdot v = A\,v\,B^*
\end{equation}
and the inner product is~$v\cdot v = \det(v)$.  

There is an equivariant 
`spinor squaring' mapping~$\sigma:\bbS^{2,2}\to\bbR^{2,2}$ defined by
$\sigma(s_+,s_-) = s_+\,s^*_-$.   
Its image is the null cone in~$\bbR^{2,2}$.

\subsubsection{Dimension 5}\label{sssec:spin5}
Again, there are three cases.

\paragraph{$\Spin(5)\simeq\Symp(2)$}\label{par:spin50}
The action of~$\Spin(5)$ on~$\bbS^{5}\simeq\bbH^2$ is
\begin{equation}
A\cdot s = A\,s.
\end{equation}
The orbits are given by the level sets of~$s\cdot s = s^*\,s$.
Except for $s=0$, these orbits all have the same stabilizer type, 
namely~$\Symp(1)$.

Identify~$\bbR^5$ with the space of traceless,
quaternion Hermitian symmetric $2$-by-$2$ matrices.  
Then the action of~$\Spin(5)$ on~$\bbR^5$ becomes
\begin{equation}
A\cdot m = A m A^*,
\end{equation}
and the quadratic form is just $m\cdot m = \tr(m^* m)$. 

There is an equivariant `spinor squaring' 
mapping~$\sigma:\bbS^{5}\to\bbR^{5}$ defined by
$\sigma(s) = s\,s^*$.  Its image is all of~$\bbR^5$.

\paragraph{$\Spin(4,1)\simeq\Symp(1,1)$}\label{par:spin41}
Let~$Q = \begin{pmatrix} 1 & 0\\ 0 & -1\end{pmatrix}$, so that~$\Spin(4,1)$
is realized as the matrices~$A\in\bbH(2)$ that satisfy~$A^*QA = Q$.
Here,~$\bbS^{4,1}\simeq\bbH^2$ and the spinor action is matrix
multiplication:
\begin{equation}
A \cdot s = A\,s\,.
\end{equation}
The spinor orbits are essentially the level sets of the 
function~$\nu:\bbS^{4,1}\to\bbR$ defined by~$\nu(s) = s^*Qs$, with the one
exception being the level set~$\nu = 0$, which consists of two orbits,
the zero vector and then everything else.  The stabilizer of 
\begin{equation}
s_0 = \begin{pmatrix} 1 \\ 1\end{pmatrix}
\qquad\text{is}\qquad
G_0 = \left\{\begin{pmatrix} 1{+}q & -q\\ -\bar q & 1+\bar q\end{pmatrix}
                   \ \vrule\ q\in\Im\bbH\ \right\}\simeq\bbR^3,
\end{equation}
while, for~$r>0$, the stabilizer of
\begin{equation}
s_{r^2} = \begin{pmatrix} r \\ 0\end{pmatrix}
 \qquad\text{is}\qquad
G_+ = \left\{\begin{pmatrix} 1 & 0\\ 0 & q\end{pmatrix}
                   \ \vrule\ q\in\Symp(1)\ \right\}\simeq \Symp(1),
\end{equation}
and the stabilizer of
\begin{equation}
s_{-r^2} = \begin{pmatrix} 0\\ r\end{pmatrix}
 \qquad\text{is}\qquad
G_- = \left\{\begin{pmatrix} q & 0\\ 0 & 1\end{pmatrix}
                   \ \vrule\ q\in\Symp(1)\ \right\}\simeq \Symp(1).
\end{equation}
The two elements~$s_{\pm r^2}$ are on the same~$\Pin(4,1)$-orbit, 
so for our purposes, they should be counted as the same.

Identify $\bbR^{4,1}$ with the space of quaternion 
Hermitian symmetric matrices~$m$ that satisfy~$\tr(Qm)=0$. Then the
action of~$\Spin(4,1)$ on this space is just
\begin{equation}
A\cdot m = A m A^*.
\end{equation}
The invariant quadratic form is~$m\cdot m = -\det(m)$, where, $\det$
is defined on the quaternion Hermitian symmetric $2$-by-$2$ matrices by
\begin{equation}
\det\begin{pmatrix}a&b\\ \bar b & c\end{pmatrix} = ac-b\bar b,
\qquad a,c\in\bbR,\ b\in\bbH.
\end{equation}

There is an equivariant 
`spinor squaring' mapping~$\sigma:\bbS^{4,1}\to\bbR^{4,1}$ defined by
$\sigma(s) = s\,s^*-\frac12\nu(s)\,Q$.   
Its image consists of half of the cone of elements~$m$ 
that satisfy~$\det(m)\ge0$.  The image boundary, i.e., 
the `forward light cone' is the image of the locus~$\nu=0$ in~$\bbS^{4,1}$.

\paragraph{$\Spin(3,2)\simeq\Symp(2,\bbR)$}\label{par:spin32}
This classical isomorphism can be described as follows:  
Let~$J = \begin{pmatrix} 0 & -\I_2 \\ \I_2 & 0\end{pmatrix}$.
Then~$\Symp(2,\bbR)$ is the subgroup of~$\GL(4,\bbR)$ consisting
of those matrices~$A$ that satisfy~$A^*JA = J$.  This group is
isomorphic to~$\Spin(3,2)$ in such a way that $\bbS^{3,2}$
can be identified with~$\bbR^{4}$ so that the spinor representation
becomes the usual matrix multiplication:
\begin{equation}
A\cdot s = As.
\end{equation}
There are only two $\Symp(2,\bbR)$-orbits in this case:  The
zero orbit and everything else.  

The vector representation is described as follows:  Identify~$\bbR^{3,2}$
with the space of skew-symmetric $v\in \bbR(4)$ that satisfy~$\tr(vJ)=0$.
This space is preserved under the action~$A\cdot v = A\,v\,A^*$.
The inner product is~$v\cdot v = \Pf(v)$.  This is an irreducible 
representation and the inner product is seen to be of type~$(3,2)$.

\subsubsection{Dimension 6}\label{sssec:spin6}
Now, there are four cases.

\paragraph{$\Spin(6)\simeq\SU(4)$}\label{par:spin60}
The action of~$\Spin(6)$ on~$\bbS^{6}\simeq\bbC^4$ is
\begin{equation}
A\cdot s = A\,s.
\end{equation}
The orbits are given by the level sets of~$s\cdot s = s^*\,s$.
Except for $s=0$, these orbits all have the same stabilizer type, 
namely~$\SU(3)$.  

To see the representation of~$\SU(4)$ on~$\bbR^{6,0}$, consider
the space~$W$ of skewsymmetric~$w\in \bbC(4)$.  This is a complex
vector space of dimension~$6$.  The group~$\SL(4,\bbC)$ acts on~$W$ 
by the rule
\begin{equation}
A\cdot w = A\,w\,A^*.
\end{equation}
Consider the complex inner product~$(,)$ on~$W$ that 
satisfies~$(w,w) = \Pf(w)$.  This is a nondegenerate quadratic form 
that is invariant under~$\SL(4,\bbC)$ and hence under~$\SU(4)$.  
There is also an Hermitian inner product on~$W$ defined by
$\la w,w\ra = \frac14\tr(ww^*)$ and it is easily seen to
be invariant under~$\SU(4)$ as well.  It follows that there is
an~$\SU(4)$-invariant conjugate-linear map~$c:W\to W$ so that
$(cw,v) = \la w, v\ra$.  This linear map satisfies~$c^2 = \I$,
so there is an $\SU(4)$-invariant splitting $W = W_+\oplus W_-$ into
the (real) eigenspaces of~$c$, each of dimension~$6$.  
The spaces~$W_\pm$ are each isomorphic 
to~$\bbR^{6,0}$ with inner product~$(,)$ and the action of ~$\SU(4)$
double covers to produce the standard~$\SO(6)$ action.

\paragraph{$\Spin(5,1)\simeq\SL(2,\bbH)$}\label{par:spin51}
Here,~$\bbS^{5,1}\simeq\bbH^2\oplus\bbH^2$ and the spinor action is
\begin{equation}
A \cdot (s_+,s_-) = (A\,s\,, (A^*)^{-1}\,s_-).
\end{equation}

There several different types of spinor orbits.
First, there is the point~$(0,0)$.   Then there are the two
orbits of dimension~$7$ of the points $(s_+,0)$ and
$(0,s_-)$, where~$s_\pm$ are nonzero.  Third, there are the orbits
that lie in the locus~$s_-^*\,s_+=0$, but that have~$s_\pm\not=0$.
These orbits all have dimension~$11$ and there is a $1$-parameter
family of them.  In fact, for each positive real~$r$, the orbit
of
\begin{equation}
s_r = \left(\begin{pmatrix}0\\1\end{pmatrix},
      \begin{pmatrix}r\\0\end{pmatrix}\right)
\qquad\text{has stabilizer}\qquad
G_0 = \left\{\ \begin{pmatrix}1&q\\0&1\end{pmatrix}\ \vrule\ 
                q\in\bbH\ \right\}\simeq\bbH.
\end{equation}
Fourth, the remaining orbits have dimension~$12$.  These are
parametrized by~$s_-^*\,s_+=\lambda\in\bbH^*$.   This level set
is the orbit of the element
\begin{equation}
s_\lambda = \left(\begin{pmatrix}1\\0\end{pmatrix},
      \begin{pmatrix}\lambda\\0\end{pmatrix}\right)
\qquad\text{with stabilizer}\qquad
G_1 = \left\{\ \begin{pmatrix}1&0\\0&q\end{pmatrix}\ \vrule\ 
                q\in\Symp(1)\ \right\}\simeq\Symp(1).
\end{equation}
Note that, because the centralizer of~$\Spin(5,1)$ 
in~$\Aut(\bbS^{5,1})$ is~$\bbH^*\times\bbH^*$ (scalar multiplication 
(on the \emph{right}) in each summand), the combined action 
of the centralizer and~$\Spin(5,1)$ shows that all of the orbits
of the third type should be regarded as essentially the same
and that all of the orbits of the fourth type should be regarded
as essentially the same.  Thus, there are really only four distinct
types of orbits to consider.  Moreover, action by an element of~$\Pin(5,1)$
not in~$\Spin(5,1)$ exchanges the two $7$-dimensional orbits, so they
should really be regarded as belonging to the same type.

Identify~$\bbR^{5,1}$ with the space of 
Hermitian symmetric $2$-by-$2$ matrices with quaternion entries.
The action of~$\Spin(5,1)$ on this space can be be described as
\begin{equation}
A\cdot a = A\,\,a\,\,A^*
\end{equation}
and the inner product satisfies~$a\cdot a = -\det(a)$, where
the interpretation of determinant in this case is 
\begin{equation}
\det\begin{pmatrix} a & b\\ \bar b & c\end{pmatrix} = ac-b\bar b
\end{equation}
for~$a,c\in\bbR$ and~$b\in\bbH$.
(That $\SL(2,\bbH)$ does preserve this must be checked, since, normally,
$\det$ is not defined for matrices with quaternion entries.)

There is an equivariant 
`spinor squaring' mapping~$\sigma_+:\bbS^{5,1}_+\to\bbR^{5,1}$ defined by
$\sigma_+(s_+) = s_+\,s^*_+$.   
Its image consists of the `forward light cone' in~$\bbR^{5,1}$.

\paragraph{$\Spin(4,2)\simeq\SU(2,2)$}\label{par:spin42}
The identification of~$\Spin(4,2)$ with~$\SU(2,2)$ is very similar with
the identification of~$\Spin(6)$ with~$\SU(4)$ and can be seen as 
follows.  Let~$Q = \begin{pmatrix}\I_2&0\\0&-I_2\end{pmatrix}$ and
recall that~$\SU(2,2)$ is the group of matrices~$A\in\SL(4,\bbC)$
satisfying~$A^*QA = Q$.  It acts on~$\bbC^{2,2} = \bbC^4$ preserving the 
Hermitian inner product defined by~$\la v,w\ra = v^*Qw$.
The orbits of this action are $0\in\bbC^4$ and the nonzero parts
of the level sets of the Hermitian form~$\la v,w\ra = v^*Qv$.
Note that the stabilizer of a vector satisfying~$v^*Qv = 0$ is not
conjugate to the stabilizer of a vector satisfying~$v^*Qv\not=0$.  
Thus, it makes sense to say that there are essentially two distinct
types of nonzero orbits, the null orbit and the non-null orbits (which
form a single type).

To justify the identification of~$\Spin(4,2)$ with~$\SU(2,2)$, it 
will be necessary to construct a $6$-dimensional real vector space~$V$
on which~$\SU(2,2)$ acts as the identity component of the stabilizer
of a quadratic form on~$V$ of type~$(4,2)$.  Here is how this can
be done:  Again, start with~$W$ being the space of skewsymmetric matrices
$w\in\bbC(4)$, with the action of~$\SL(4,\bbC)$ being, as before,
$A\cdot w = A\,w\,A^*$. Again define the complex inner product~$(,)$ on~$W$ 
so that~$(w,w) = \Pf(w)$.  Now, consider the Hermitian inner product on~$W$
defined by~$\la w, v\ra = \frac14\tr(w^*Qv)$. 
This Hermitian inner product is invariant under~$\SU(2,2)$, so there
is an $\SU(2,2)$-invariant conjugate linear mapping~$c:W\to W$ 
satisfying~$(cw,v)=\la w,v\ra$.
Again,~$c^2$ is the identity, so that~$W$ can be split into real subspaces
$W = W_+\oplus W_-$ with~$i\,W_\pm = W_\mp$.  Then~$\SU(2,2)$ acts on
$V = W_+$ preserving~$(,)$ and it is not difficult to see that the type of
this quadratic form is~$(4,2)$.  Since~$\SU(2,2)$ is simple and of 
dimension~$15$, the same dimension as~$\SO(4,2)$, it follows that this
representation of~$\SU(2,2)$ must be onto the identity component of
the stabilizer of this quadratic form, as desired.  More detail about 
this representation will be supplied when it is needed in the next 
section.

\paragraph{$\Spin(3,3)\simeq\SL(4,\bbR)$}\label{par:spin33}
Here~$\bbS^{3,3}\simeq\bbR^4\oplus\bbR^4$ and the spinor action is
\begin{equation}
A \cdot (s_+,s_-) = (As_+,\,(A^*)^{-1}s_-).
\end{equation}
There are several orbits of~$\Spin(3,3)$ on~$\bbS^{3,3}$:
Those of the points~$(0,0)$, $(s_+,0)$, 
$(0,s_-)$, and $(s_+,s_-)$ where $s^*_-s_+=\lambda$, where~$\lambda$ is
any real number and~$s_\pm$ are nonzero elements of~$\bbR^4$.  
In this last family of orbits, there are two essentially different kinds.  
The orbit with~$\lambda=0$ has a different stabilizer type in~$\SL(4,\bbR)$ 
from those with~$\lambda\not=0$, even though it has the same dimension.  
This is accounted for by the fact that the centralizer of~$\Spin(3,3)$ 
in~$\Aut(\bbS^{3,3})$ is~$\bbR^*\times\bbR^*$ (scalar multiplication in
the fibers) and the combined action of the centralizer and~$\Spin(3,3)$
makes all of the orbits with $\lambda\not=0$ equivalent to each other.
Moreover, action by an element of~$\Pin(3,3)$ not in~$\Spin(3,3)$
exchanges the two $4$-dimensional orbits, so they should be 
regarded as belonging to the same orbit type.

Under the identification~$\bbR^{3,3}\simeq A_4(\bbR)$, the 
antisymmetric $4$-by-$4$ matrices with real entries, the action of~$\Spin(3,3)$
can be be described as
\begin{equation}
A\cdot a = A\,\,a\,\,A^*
\end{equation}
and the inner product satisfies~$a\cdot a = \Pf(a)$.

\subsection{The split cases and pure spinors}\label{ssec:spin_pure}
The orbit structure of~$\Spin(p,q)$ grows increasingly
complicated as~$p{+}q$ increases.  However, there are a few
orbits that are easy to describe in the so-called `split' cases,
i.e., $\Spin(p{+}1,p)$ (the odd split case), and~$\Spin(p,p)$ (the
even split case).  

When $p = q$ or~$p = q{+}1$, the maximal dimension 
of a null plane~$N\subset\bbR^{p,q}$ is~$q$.  Let~$v_1,\ldots,v_q$
be a basis of such an~$N$ and let~$[v] = v_1v_2\cdots v_q\in\Cl(p,q)$.
The element~$[v]$ depends only on~$N$ and a choice of volume element for~$N$.
It is not hard to show that the left ideal~$\Cl(p,q)\cdot[v]\subset\Cl(p,q)$ 
is minimal, and so is irreducible as a~$\Cl(p,q)$ module.

\subsubsection{The odd case}\label{sssec:spin_pure_odd}
Now, according to the definitions in~\S\ref{ssec:spinpq},
when~$p = q{+}1$, the odd case, $\Cl(q{+}1,q)\cdot[v]$,
when considered as a $\Spin(q{+}1,q)$-module, is two
isomorphic copies of~$\bbS^{q+1,q}$.  Fix such a decomposition
of~$\Cl(p,q)\cdot[v]$ and consider the image~$\la v\ra$
of~$[v]$ in one of these summands, henceforth denoted~$\bbS^{q+1,q}$.
The $\Spin(q{+}1,q)$-orbit of~$\la v\ra$ is known
as the space of \emph{pure spinors}.  This orbit is a cone and
has dimension~$\frac12 q(q{+}1) + 1$, which turns out to be the lowest 
dimension possible for a nonzero orbit.  The $\rho$-image of
the stabilizer in~$\Spin(q{+}1,q)$ of a pure spinor 
is the stabilizer in~$\SO(q{+}1,q)$ 
of a corresponding null $q$-vector in~$\bbR^{p,p}$.  

\paragraph{Low values of~$q$}\label{par:spin_pure_odd_low}
When~$q\equiv 0,3\mod 4$, $\Spin(q{+}1,q)$ preserves an inner product
(of split type) on~$\bbS^{q+1,q}$ while, when~$q\equiv 1,2\mod 4$, 
$\Spin(q{+}1,q)$ preserves a symplectic form on~$\bbS^{q+1,q}$, 
see~\cite{rHa90}.

Since~$\bbS^{q+1,q}$ is a real vector space of dimension~$2^q$, 
as~$q$ increases, the pure spinors become a relatively small
$\Spin(q{+}1,q)$-orbit in~$\bbS^{q+1,q}$.  

However, for low values of~$q$, the situation is different.  
When~$q=1$~or~$2$, every spinor is pure.  

When~$q=3$, dimension count shows that the pure spinors are a hypersurface 
in~$\bbS^{4,3}$.  Since they form a cone, they must constitute 
the null cone in~$\bbS^{4,3}\simeq\bbR^8$ of the~$\Spin(4,3)$-invariant 
quadratic form on~$\bbS^{4,3}$.  
Moreover, the other nonzero $\Spin(4,3)$-orbits 
in~$\bbS^{4,3}$ are the nonzero level sets of this quadratic form,
and so are also of dimension~$7$.  The stabilizer of a non-null
element~$v\in\bbS^{4,3}$ is isomorphic to~$\G^*_2\subset\Spin(4,3)$,
the split form of type~$\G_2$.  

When~$q=4$, the pure spinors constitute an $11$-dimensional cone
in~$\bbS^{5,4}\simeq\bbR^{16}$, which must therefore lie in the null 
cone of the $\Spin(5,4)$-invariant quadratic form on~$\bbS^{5,4}$.  
It is an interesting fact that each of the nonzero level sets of this 
quadratic form constitutes a single~$\Spin(5,4)$-orbit.  (This is
because, as can be seen in~\cite{rBr99spinors}, 
$\Spin(9)$ acts transitively on
the unit spheres in~$\bbS^{9}\simeq\bbR^{16}$.  The existence of
hypersurface orbits in the compact case implies the existence of 
hypersurface orbits in the complexification, which implies the existence 
of hypersurface orbits in the split form, i.e., $\Spin(5,4)$.)
Thus, although the null cone is the limit of hypersurface orbits,
it does not constitute a single orbit, but must contain at least two
orbits (besides the zero orbit).  One of those orbits is the $11$-dimensional
space of pure spinors, but I do not know whether the complement of
the pure spinors in the null spinors constitutes a single orbit or not.


\subsubsection{The even case}\label{sssec:spin_pure_even}
According to the definitions in~\S\ref{ssec:spinpq}, when~$p=q$,
the relation $\bbS^{p,p}_+\oplus \bbS^{p,p}_-= \bbS^{p,p} \simeq
\Cl(p,p)\cdot[v]$ holds.  It turns out that $[v]$ lies in one of the
two summands (which one depends on the orientation of~$\bbR^{p,p}$,
since this decides which one is~$\bbS^{p,p}_+$).  This corresponds to
the well-known fact that the space of maximal null $p$-planes 
in~$\bbR^{p,p}$ consists of two components.  By this construction, 
each component of the space of null $p$-planes endowed with a 
choice of volume form in~$\bbR^{p,p}$ is double covered by a $\Spin(p,p)$
orbit (in fact, a closed cone) in~$\bbS^{p,p}_\pm$.  The elements of
these two orbits are the pure spinors.  Each forms a
minimal (i.e., maximally degenerate) orbit in ~$\bbS^{p,p}$.  The
dimension of each of these orbits is $\frac12 p(p{-}1) + 1$. 
The $\rho$-image of the stabilizer in~$\Spin(p,p)$ of a pure spinor 
maps onto the stabilizer of a null $p$-vector in~$\bbR^{p,p}$.  

\paragraph{Low values of~$p$}\label{par:spin_pure_even_low}
When~$p\equiv 1\mod 2$, the spaces~$\bbS^{p,p}_+$ and~$\bbS^{p,p}_-$
are naturally dual as $\Spin(p,p)$-modules.  When~$p\equiv 2\mod 4$, 
each of~$\bbS^{p,p}_\pm$ is a symplectic representation of~$\Spin(p,p)$.
When~$p\equiv 0\mod 4$, each of~$\bbS^{p,p}_\pm$ is an orthogonal
representation of~$\Spin(p,p)$. Again, see~\cite{rHa90} for proofs
of these facts.

Since~$\bbS^{p,p}$ is a sum of two $\Spin(p,p)$-irreducible 
real vector spaces of dimension~$2^{p-1}$, 
as~$p$ increases, the pure spinors become a vanishingly small
$\Spin(p,p)$-orbit in~$\bbS^{p,p}$.  

However, for low values of~$p$, the situation is different.  
When~$p=1$,~$2$, or~$3$, every spinor in~$\bbS^{p,p}_\pm$ is pure.  

When~$p=4$ (the famous case of triality), $\Spin(4,4)$ acts on
each of~$\bbS^{4,4}_\pm\simeq\bbR^{4,4}$ as the full group of linear 
transformations preserving the spinor inner product.  In particular, 
the nonzero orbits are just the level sets of the invariant quadratic 
form. Thus, the pure spinors in each space constitute the
null cone (minus the origin) of the quadratic form.  Using this
description, it is not difficult completely to describe the orbits
of~$\Spin(4,4)$ on~$\bbS^{4,4}$.  I will go into more detail 
as necessary in what follows.

When~$p=5$, the situation is more subtle.  $\Spin(5,5)$ acts
on each of~$\bbS^{5,5}_\pm\simeq\bbR^{16}$ with open orbits.
The cone of pure spinors in each summand has dimension~$11$.  
In fact, in the direct sum action on~$\bbS^{5,5}$, the
group $\Spin(5,5)$ preserves the quadratic form that is the 
dual pairing on the two factors and a nontrivial quartic form.
The generic orbits of ~$\Spin(5,5)$ on~$\bbS^{5,5}$ are simultaneous
level sets of these two polynomials and so have dimension~$30$.
I do not know the full orbit structure.

\subsection{The octonions and $\Spin(10,1)$}\label{ssec:spin_oct+spin10,1}
In this section, I will develop just enough of the necessary algebra
to discuss the geometry of one higher dimensional case, that of
parallel spinors in a metric of type~$(10,1)$.  The reason 
for considering this case is that there is some interest in it
for physical reasons, see~\cite{jFF99}.

\subsubsection{Octonions}\label{sssec:octonions}
A few background facts about the octonions will be needed.  For
proofs, see~\cite{rHa90}.

As usual, let~$\bbO$ denote the ring of octonions. 
Elements of~$\bbO$ will be denoted by bold letters, such as $\xb$, $\yb$, 
etc.  Thus, $\bbO$ is the unique $\bbR$-algebra of dimension~$8$
with unit~$\oneb\in\bbO$ endowed with a positive definite inner 
product~$\la,\ra$ satisfying~$\la \xb\yb,\xb\yb\ra
=\la\xb,\xb\ra\,\la\yb,\yb\ra$ for all~$\xb,\yb\in\bbO$.   
As usual, the norm of
an element~$\xb\in\bbO$ is denoted $|\xb|$ and defined as the square root
of~$\la\xb,\xb\ra$.  Left and right multiplication
by~$\xb\in\bbO$ define maps~$L_\xb\,,R_\xb:\bbO\to\bbO$ that are isometries
when~$|\xb|=1$.

The conjugate of~$\xb\in\bbO$, 
denoted~$\overline\xb$, is defined to be~$\overline\xb = 
2\la\xb,\oneb\ra\,\oneb-\xb$.  When a symbol is needed, the map of 
conjugation will be denoted $C:\bbO\to\bbO$.  The identity
$\xb\,\overline\xb = |\xb|^2$ holds, as well as the conjugation 
identity~$\overline{\xb\yb}={\overline\yb}\,{\overline\xb}$.
In particular, this implies the useful identities 
$C\, L_\xb\, C=R_{\overline\xb}$
and $C\, R_\xb\, C=L_{\overline\xb}$.

The algebra~$\bbO$ is not commutative or associative.  However, any 
subalgebra of~$\bbO$ that is generated by two elements is associative.  
It follows that $\xb\,\bigl({\overline\xb} \yb\bigr)=|\xb|^2\,\yb$ and
that~$(\xb\yb)\xb=\xb(\yb\xb)$ for 
all~$\xb,\yb\in\bbO$. Thus,~$R_\xb\, L_\xb
=L_\xb\, R_\xb$ (though, of course, 
$R_\xb\, L_\yb\not=L_\yb\, R_\xb$ in general). In
particular, the expression $\xb\yb\xb$ is unambiguously defined.
In addition, there are the {\it Moufang Identities}
\begin{equation}
\begin{split}
(\xb\yb\xb)\zb &= \xb\bigl(\yb(\xb\zb)\bigr),\\
\zb(\xb\yb\xb) &= \bigl((\zb\xb)\yb\bigr)\xb,\\
\xb(\yb\zb)\xb &= (\xb\yb)(\zb\xb),
\end{split}\label{eq:MoufangIdentities}
\end{equation}
which will be useful below. 

\subsubsection{$\Spin(8)$}\label{sssec:spin8}
For~$\xb\in \bbO$, define the linear 
map~$m_\xb:\bbO\oplus\bbO\to\bbO\oplus\bbO$ by the formula
\begin{equation}
m_\xb = \left[\begin{matrix} 0&C \, R_\xb\\ 
                             -C \, L_\xb&0\end{matrix}\right]\,.
\end{equation}
By the above identities, it follows that $(m_\xb)^2 = -|\xb|^2$ and
hence this map induces a representation on the vector 
space~$\bbO\oplus\bbO$ of the Clifford algebra generated 
by~$\bbO$ with its standard quadratic form.  
This Clifford algebra is known to be isomorphic
to~$M_{16}(\bbR)$, the algebra of $16$-by-$16$ matrices with real entries, so
this representation must be faithful.  By dimension count, this establishes
the isomorphism~$\Cl\bigl(\bbO,\la,\ra\bigr)
=\End_\bbR\bigl(\bbO\oplus\bbO\bigr)$.

The group $\Spin(8)\subset\GL_\bbR(\bbO\oplus\bbO)$ is defined as the
subgroup generated by products of the form~$m_\xb\,m_\yb$
where~$\xb,\yb\in\bbO$ satisfy~$|\xb|=|\yb|=1$.  
Such endomorphisms preserve the splitting of $\bbO\oplus\bbO$ into
the two given summands since
\begin{equation}
m_\xb\,m_\yb 
= \left[\begin{matrix} -L_{\overline\xb} \, L_\yb&0\\ 
                 0&-R_{\overline\xb} \, R_\yb\end{matrix}\right]\,.
\end{equation}
In fact, setting $\xb=-\oneb$ in this formula shows that endomorphisms
of the form
\begin{equation}
\left[\begin{matrix} L_\ub&0\\ 0&R_\ub\end{matrix}\right],
\qquad\text{with $|\ub|=1$}
\end{equation}
lie in~$\Spin(8)$.  In fact, they generate~$\Spin(8)$, since 
$m_\xb\,m_\yb$ is clearly a product of two of these when 
$|\xb|=|\yb|=1$.

Fixing an identification~$\bbO\simeq\bbR^8$ defines an embedding
$\Spin(8)\subset\SO(8)\times\SO(8)$, and the projections onto either
of the factors is a group homomorphism.  Since neither of these
projections is trivial, since the Lie algebra~$\euso(8)$ is simple, and
since~$\SO(8)$ is connected, it follows that each of these projections
is a surjective homomorphism.  Since~$\Spin(8)$ is simply connected
and since the fundamental group of~$\SO(8)$ is~$\bbZ_2$, it follows that  
that each of these homomorphisms is a non-trivial double cover of~$\SO(8)$.
Moreover, it follows that the subsets~$\{\  L_\ub\ \vrule\ |\ub|=1\ \}$
and $\{\  R_\ub\ \vrule\ |\ub|=1\ \}$ of~$\SO(8)$ each suffice to
generate~$\SO(8)$.

Let~$H\subset\bigl(\SO(8)\bigr)^3$ be the set of 
triples~$(g_1,g_2,g_3)\in\bigl(\SO(8)\bigr)^3$ for which
\begin{equation}
g_2(\xb\yb) = g_1(\xb)\,g_3(\yb)
\end{equation}
for all~$\xb,\yb\in\bbO$.  The set~$H$ is closed and is evidently 
closed under multiplication and inverse. Hence it is a compact Lie group.  

By the third Moufang identity, $H$ contains the subset
\begin{equation}
\Sigma = 
\left\{\ (L_\ub,\,L_\ub{\,}R_\ub,\,R_\ub)\ \vrule\ |\ub|=1
\right\}.  
\end{equation}
Let~$K\subset H$ be the subgroup generated by~$\Sigma$, and for $i=1,2,3$,
let~$\rho_i:H\to\SO(8)$ be the homomorphism that is 
projection onto the $i$-th factor.  Since~$\rho_1(K)$ contains
$\{\  L_\ub\ \vrule\ |\ub|=1\ \}$, it follows that $\rho_1(K)=\SO(8)$,
so, \emph{a fortiori}, $\rho_1(H)=\SO(8)$.  Similarly, $\rho_3(H)=\SO(8)$.  

The kernel of~$\rho_1$ consists of elements~$(I_8,g_2,g_3)$
that satisfy~$g_2(\xb\yb) = \xb\,g_3(\yb)$ for all~$\xb,\yb\in\bbO$.  
Setting~$\xb=\oneb$ in this equation yields~$g_2=g_3$,
so that~$g_2(\xb\yb) = \xb\,g_2(\yb)$.  Setting~$\yb=\oneb$
in this equation yields~$g_2(\xb) = \xb\,g_2(\oneb)$, i.e., $g_2 = R_\ub$
for $\ub=g_2(\oneb)$.  Thus, the elements in the kernel of~$\rho_1$ are
of the form~$(1,R_\ub,R_\ub)$ for some~$\ub$ with~$|\ub|=1$.  
However, any such $\ub$ would, by definition, 
satisfy~$(\xb\yb)\ub=\xb(\yb\ub)$ for all~$\xb,\yb\in\bbO$, 
which is impossible unless~$\ub=\pm\oneb$.  Thus, the kernel of~$\rho_1$
is~$\bigl\{(I_8,\pm I_8,\pm I_8)\bigr\}\simeq\bbZ_2$, so that $\rho_1$ is a 
$2$-to-$1$ homomorphism of~$H$ onto~$\SO(8)$.  Similarly, $\rho_3$ is a 
$2$-to-$1$ homomorphism of~$H$ onto~$\SO(8)$, with 
kernel~$\bigl\{(\pm I_8,\pm I_8,I_8)\bigr\}$.  Thus, $H$ is either 
connected and isomorphic to~$\Spin(8)$ or else disconnected,
with two components.

Now~$K$ is a connected subgroup of~$H$ and the kernel of~$\rho_1$ 
intersected with~$K$ is either trivial or~$\bbZ_2$.  Moreover, the product 
homomorphism~$\rho_1{\times}\rho_3:K\to\SO(8){\times}\SO(8)$ maps the 
generator~$\Sigma\subset K$ into generators 
of~$\Spin(8)\subset\SO(8){\times}\SO(8)$.  It follows that
$\rho_1{\times}\rho_3(K)=\Spin(8)$ and hence that~$\rho_1$ and $\rho_3$ 
must be non-trivial double covers of~$\Spin(8)$ when restricted to~$K$.  
In particular, it follows that~$K$ must be all of~$H$ and, moreover, that the 
homomorphism~$\rho_1{\times}\rho_3:H\to\Spin(8)$ must be an isomorphism.
It also follows that the homomorphism $\rho_2:H\to\SO(8)$ 
must be a double cover of~$\SO(8)$ as well.  

Henceforth, $H$ will be identified with~$\Spin(8)$ via the 
isomorphism~$\rho_1{\times}\rho_3$.  Note that the center of $H$
consists of the 
elements~$(\varepsilon_1\,I_8,\varepsilon_2\,I_8,\varepsilon_3\,I_8)$
where~${\varepsilon_i}^2=\varepsilon_1\varepsilon_2\varepsilon_3=1$
and is isomorphic to~$\bbZ_2\times\bbZ_2$.

\paragraph{Triality}\label{par:triality}
For~$(g_1,g_2,g_3)\in H$, the identity
$g_2(\xb\yb) = g_1(\xb)\,g_3(\yb)$ can be conjugated, giving
\begin{equation}
Cg_2C(\xb\yb) 
= \overline{g_2(\overline{\yb}\,\overline{\xb})}
= \overline{g_1(\overline{\yb})\,g_3(\overline{\xb})}
= \overline{g_3(\overline{\xb})}\, \overline{g_1(\overline{\yb})}.
\end{equation}
This implies that $\bigl(Cg_3C,Cg_2C,Cg_1C\bigr)$ also lies in~$H$.
Also, replacing $\xb$ by $\zb\overline{\yb}$ in the original formula
and multiplying on the right by $\overline{g_3(\yb)}$ shows that
\begin{equation}
g_2(\zb) \overline{g_3(\yb)} = g_1(\zb\overline{\yb}),
\end{equation}
implying that~$\bigl(g_2,g_1,Cg_3C\bigr)$ lies in~$H$ as well.  In fact, the
two maps~$\alpha,\beta:H\to H$ defined by
\begin{equation}
\alpha(g_1,g_2,g_3) = \bigl(Cg_3C,Cg_2C,Cg_1C\bigr),
\quad\hbox{and}\qquad
\beta(g_1,g_2,g_3) = \bigl(g_2,g_1,Cg_3C\bigr) 
\end{equation}
are outer automorphisms (since they act nontrivially on the center of~$H$) 
and generate a group of automorphisms isomorphic to~$S_3$, the symmetric group 
on three letters.  The automorphism $\tau=\alpha\beta$ is known as the 
triality automorphism.

To emphasize the group action, denote~$\bbO\simeq\bbR^8$ by~$V_i$
when regarding it as a representation space of~$\Spin(8)$ via
the representation~$\rho_i$.  Thus, octonion multiplication induces
a $\Spin(8)$-equivariant projection
\begin{equation}
V_1\otimes V_3 \longrightarrow V_2\,.
\end{equation}
In the standard notation, it is traditional to identify $V_1$ with
$\bbS^8_-$ and~$V_3$ with $\bbS^8_+$ and to refer to $V_2$ as the 
`vector representation'~$\bbR^8$.  Let~$\rho_i':\euspin(8)\to\euso(8)$ 
denote the corresponding Lie algebra homomorphisms, which are, in fact, 
isomorphisms.  For simplicity of notation, for any~$a\in\euspin(8)$,
the element~$\rho_i'(a)\in\euso(8)$ will be denoted by~$a_i$ when
no confusion can arise.

\subsubsection{$\Spin(10,1)$}\label{sssec:spin101}
I will now go directly to the construction of~$\Spin(10,1)$ and its
usual spinor representation.  For more detail and for justification
of some of the statements, the reader can consult~\cite{rBr99spinors},
although there are, of course, many classical sources for this material.

It is convenient to identify~$\bbC\otimes\bbO^2$
with~$\bbO^4$ explicitly via the identification
\begin{equation}
\zb = \begin{pmatrix}\xb_1+i\,\xb_2\\ \yb_1+i\,\yb_2\end{pmatrix}
    = \begin{pmatrix}\xb_1\\ \yb_1\\ \xb_2\\ \yb_2\end{pmatrix}.
\end{equation}
Via this identification, $\euspin(10)$ can be identified with the subspace
\begin{equation}
\euspin(10) 
= \left\{ \begin{pmatrix}
   a_1    & C\,R_\xb &  -r\,I_8  &  -C\,R_\yb\\
-C\,L_\xb & a_3      & -C\,L_\yb &  r\,I_8   \\
r\,I_8    & C\,R_\yb &    a_1    &   C\,R_\xb\\
 C\,L_\yb &  -r\,I_8 & -C\,L_\xb &     a_3   \end{pmatrix}\ 
  \vrule\ \begin{matrix}r\in\bbR,\\ \xb,\yb\in\bbO,\\ a\in\euspin(8)\end{matrix}
   \right\}\,.
\end{equation}
Consider the one-parameter subgroup~$\bR\subset\SL_\bbR(\bbO^4)$ defined by
\begin{equation}
\bR 
=\left\{\begin{pmatrix}t\,I_{16}&0\\0&t^{-1}\,I_{16}\end{pmatrix}
  \ \vrule\ t\in\bbR^+\ \right\}.
\end{equation}
It has a Lie algebra~$\eur\subset\eusl(\bbO^4)$. Evidently, the
the subspace~$\left[\euspin(10),\eur\right]$ consists of matrices of
the form
\begin{equation}
\begin{pmatrix}
   0_8 & 0_8 &  r\,I_8  &  C\,R_\yb\\
   0_8 & 0_8 & C\,L_\yb &  -r\,I_8   \\ 
 r\,I_8    & C\,R_\yb & 0_8 & 0_8 \\
 C\,L_\yb &  -r\,I_8 & 0_8 &  0_8 \end{pmatrix},
\qquad r\in\bbR,\ \yb\in\bbO\,.
\end{equation}
Let~$\eug = \euspin(10) \oplus\eur\oplus\left[\euspin(10),\eur\right]$. 
Explicitly,
\begin{equation}
\eug 
= \left\{\begin{pmatrix}
   a_1 +x\,I_8  & C\,R_\xb &  y\,I_8  &  C\,R_\yb\\
-C\,L_\xb & a_3 +x\,I_8  & C\,L_\yb &  -y\,I_8   \\ 
z\,I_8    & C\,R_\zb &  a_1 -x\,I_8   &   C\,R_\xb\\
 C\,L_\zb &  -z\,I_8 & -C\,L_\xb &     a_3 -x\,I_8   \end{pmatrix}
\ \vrule\ 
\begin{matrix}x,y,z\in\bbR,\\ \noalign{\vskip2pt} \xb,\yb,\zb\in\bbO,\\
        \noalign{\vskip2pt}\ a\in\euspin(8)\end{matrix}\ \right\}\,.
\end{equation}
One can show that $\eug$ is isomorphic to~$\euso(10,1)$ 
and hence is the Lie algebra of a representation 
of~$\Spin(10,1)$.  It is not hard to argue that this representation
on~$\bbO^4\simeq\bbR^{32}$ must be equivalent to the 
representation~$\bbS^{10,1}$.

Thus, define~$\Spin(10,1)$ to be the (connected) subgroup 
of~$\SL_\bbR(\bbO^4)$ that is generated by~$\Spin(10)$ and the 
subgroup~$\bR$.  Its Lie algebra~$\eug$ will henceforth be written 
as~$\euspin(10,1)$.

Consider the polynomial
\begin{equation}
p(\zb) = |\xb_1|^2|\xb_2|^2+|\yb_1|^2|\yb_2|^2 
     - \left(\xb_1\cdot\xb_2+\yb_1\cdot\yb_2\right)^2
     + 2\,(\xb_1\yb_1)\cdot(\xb_2\yb_2)\,.
\end{equation}
It is not difficult to show that $p$ is nonnegative and is also 
invariant under the action of~$\Spin(10,1)$. Moreover, the orbits
of~$\Spin(10,1)$ are the positive level sets of this polynomial
and the zero level set minus the origin.  The positive level sets
are smooth and have dimension~$31$, while the zero level set is
smooth away from the origin and has dimension~$25$.  

In fact, $p$ has the following interpretation:
Consider the squaring map~$\sigma:\bbO^4\to\bbR^{2,1}\oplus\bbO=\bbR^{10,1}$ 
that takes spinors for~$\Spin(10,1)$ to vectors.  
This map~$\sigma$ is defined as follows:
\begin{equation}
\sigma\left(\begin{pmatrix}\xb_1\\ \yb_1\\ \xb_2\\ \yb_2\end{pmatrix}\right) 
= \begin{pmatrix}|\xb_1|^2+|\yb_1|^2\\ 
        2\,\bigl(\xb_1\cdot\xb_2-\yb_1\cdot\yb_2\bigr)\\
        |\xb_2|^2+|\yb_2|^2\\
      2\,\bigl(\xb_1\,\yb_2+\xb_2\,\yb_1\bigr)\end{pmatrix}.
\end{equation}
Define the inner product on vectors in~$\bbR^{2,1}\oplus\bbO=\bbR^{10,1}$ 
by the rule
\begin{equation}
\begin{pmatrix}a_1\\ a_2\\ a_3\\ \xb\end{pmatrix}\cdot
\begin{pmatrix}b_1\\ b_2\\ b_3\\ \yb\end{pmatrix} 
= -2(a_1b_3+a_3b_1)+a_2b_2 + \xb\cdot\yb
\end{equation}
and let~$\SO(10,1)$ denote the subgroup of~$\SL(\bbR^{2,1}\oplus\bbO)$
that preserves this inner product.  This group still has two components
of course, but only the identity component~$\SO^\uparrow(10,1)$ will be of 
interest here. Let~$\rho:\Spin(10,1)\to\SO^\uparrow(10,1)$ be the homomorphism 
whose induced map on Lie algebras is given by the isomorphism
\begin{equation}
\rho'\left(
\begin{pmatrix}
   a_1 +x\,I_8  & C\,R_\xb &  y\,I_8  &  C\,R_\yb\\
-C\,L_\xb & a_3 +x\,I_8  & C\,L_\yb &  -y\,I_8   \\ 
z\,I_8    & C\,R_\zb &  a_1 -x\,I_8   &   C\,R_\xb\\
 C\,L_\zb &  -z\,I_8 & -C\,L_\xb &     a_3 -x\,I_8   \end{pmatrix}
\right)
= 
\begin{pmatrix}
   2x  & y &  0  &  \overline{\yb}^*\\
2z & 0  & 2y &  2\,\overline{\xb}^*   \\ 
0   & z &   -2x   &   \overline{\zb}^*\\
 2\,\overline{\zb} &  -2\,\overline{\xb} & 2\,\overline{\yb} & a_2\end{pmatrix}.
\end{equation}
The map~$\sigma$ has the
equivariance~$\sigma\bigl(g\,\zb\bigr) = \rho(g)\,\bigl(\sigma(\zb)\bigr)$ 
for ~$g\in\Spin(10,1)$ and~$\zb\in\bbO^4$.  

With these definitions, the polynomial~$p$ has the expression~$p(\zb) 
= -\frac14\,\sigma(\zb)\cdot\sigma(\zb)$, from which its invariance
is immediate.  Moreover, it follows from this that~$\sigma$ 
carries the orbits of~$\Spin(10,1)$ to the orbits of~$\SO^\uparrow(10,1)$
and that the image of~$\sigma$ is the union of the origin, the forward
light cone, and the future-directed time-like vectors.

In particular, a spinor~$\zb$ that satisfies~$p(\zb)>0$ defines a
non-zero time-like vector~$\sigma(\zb)\in\bbR^{10,1}$.  Using this fact,
it follows without difficulty that the stabilizer of such a~$\zb$ is
a conjugate of~$\SU(5)\subset\Spin(10)\subset\Spin(10,1)$.  On the other
hand, the Lie algebra~$\euh$ of the stabilizer for the null spinor
\begin{equation}
\zb_0 = \begin{pmatrix}\oneb\\0\\0\\0\end{pmatrix}
\qquad\text{is}\qquad
\euh = \left\{\begin{pmatrix}
a_1& 0 &  y\,I_8  &  C\,R_\yb\\
 0 &a_3& C\,L_\yb &  -y\,I_8   \\ 
 0 & 0 &a_1& 0 \\
 0 & 0 & 0 &a_3\end{pmatrix}\ 
  \vrule\ 
\begin{matrix}y\in\bbR,\\\noalign{\vskip2pt} \yb\in\bbO,\\
        \noalign{\vskip2pt}\ a\in\euk_1\end{matrix}\ \right\},
\end{equation}
where~$\euk_1$ is the Lie algebra of~$K_1\subset\Spin(8)$.  Thus, the
stabilizer is a semi-direct product of~$\Spin(7)$ with a copy of~$\bbR^9$,
and so has dimension~$30=55-25$, as desired.  

In conclusion, there are essentially two distinct types of $\Spin(10,1)$
orbits in~$\bbS^{10,1}$, those of the positive level sets of~$p$ and
the nonzero elements in the zero level set of~$p$.

\section{Metrics with Parallel Spinor Fields}\label{sec:parspin}

In this section, I will describe some of the normal forms and methods
for obtaining them for metrics that have parallel spinor fields.

\subsection{Dimension~$3$}\label{ssec:parspin3}

As a warmup, consider the case of metrics in dimension~$3$.

\subsubsection{Type~$(3,0)$}\label{sssec:parspin30}
Recall that~$\Spin(3)\simeq\Symp(1)$, with~$\bbS^{3,0}\simeq\bbH$. 
Thus, the $\Spin(3)$-stabilizer of any nonzero element 
of~$\bbS^{3,0}$ is trivial.  Consequently, if~$(M^3,g)$ has
a nonzero parallel spinor field, its holonomy is trivial and the
metric is flat.

\subsubsection{Type~$(2,1)$}\label{sssec:parspin21}
Since~$\Spin(2,1)$
is isomorphic to~$\SL(2,\bbR)$, with~$\bbS^{2,1}\simeq\bbR^2$, all of the
nonzero spinors constitute a single orbit.  In particular, the stabilizers
of these are all conjugate to the one-dimensional unipotent upper triangular
matrices in~$\SL(2,\bbR)$.  Thus, take the
structure equations for coframes~$\omega_{ij}=\omega_{ji}$ so that
\begin{equation}
g = \omega_{11}\,\omega_{2} 
      - \omega_{21}\,\omega_{12} 
    = \omega_{11}\,\omega_{2} - {\omega_{21}}^2
\end{equation} 
to have the form
\begin{equation}
d\begin{pmatrix}\omega_{11}&\omega_{12}\\ 
                \omega_{21}&\omega_{22} \end{pmatrix}
= -\begin{pmatrix}0&\alpha\\0&0\end{pmatrix}\w
   \begin{pmatrix}\omega_{11}&\omega_{12}\\ 
        \omega_{21}&\omega_{22}\end{pmatrix}
  +\begin{pmatrix}\omega_{11}&\omega_{12}\\ 
      \omega_{21}&\omega_{22}\end{pmatrix}
    \w \begin{pmatrix}0&0\\ \alpha&0\end{pmatrix}.
\end{equation}
Since~$d\omega_{22}=0$, I can write~$\omega_{22}=dx_{22}$
for some function~$x_{22}$.  Since $d\omega_{21}
=\omega_{22}\w\alpha$, there exists locally a  
coordinate~$x_{21}$ so that~$\omega_{21}=dx_{21}-p\,dx_{22}$.  
This makes $\alpha = dp + q\,dx_{22}$ for some function~$q$.
Reducing frames to make~$p=0$ (which can clearly be done) 
makes~$\alpha = q\,dx_{22}$ and
\begin{equation}
d\omega_{11}=-2\,\alpha\w\omega_{21}
    = 2\,q\,dx_{21}\w dx_{22}\,,
\end{equation}
so that there must be a function~$f$ on an open set in~$\bbR^2$ so that
\begin{equation}
2\,q\,dx_{21}\w dx_{22}
= d\bigl(f(x_{21},x_{22})\, dx_{22}\bigr).
\end{equation}
Thus, there is an $\bbR$-valued coordinate~$x_{11}$ so that
$\omega_{11} = dx_{11} + f(x_{21},x_{22})\, dx_{22}$.  
In particular, the metric~$g$ is locally of the form
\begin{equation}
g = dx_{11}{\circ}dx_{22} 
      - dx_{21}{\circ}dx_{12} 
   + f(x_{21},x_{22})\, (dx_{2\bar2})^2.
\end{equation}

Conversely, via this formula, any function~$f$ of two variables will 
produce a $(2,1)$-metric with a parallel spinor field.  Note that $g$
will be flat if and only if the curvature $2$-form
\begin{equation}
F = d\alpha 
= d\left(\frac12\,\frac{\p f}{\p x_{21}}\ dx_{22}\right)
\end{equation}
vanishes.   Of course, imposing the Einstein condition makes the 
curvature vanish identically.

Since the ambiguity in the choice of coordinates~$x_{22},x_{21},x_{11}$
involved only choosing arbitrary functions of one variable, it makes sense 
to say that the general metric of type~$(2,1)$ that has a parallel
spinor field depends on one function of two variables.

\subsection{Dimension~$4$}\label{ssec:parspin4}

In this subsection, I will review the well-known 
classification of pseudo-Riemannian 
metrics with parallel spinors in dimension~$4$.

\subsubsection{Type~$(4,0)$}\label{sssec:parspin40}
Since~$\Spin(4)\simeq\Symp(1)\times\Symp(1)$ and there are only
two orbit types (up to orientation), there are only two
possibilities:

\paragraph{Generic}\label{par:parspin40gen}
If $(M^4,g)$ has a parallel spinor of generic type, then
its holonomy is a subgroup of the stabilizer of the generic type, 
i.e., it is trivial, so~$(M^4,g)$ is flat.   

\paragraph{Special}\label{par:parspin40special}
If $(M^4,g)$ has a nonzero parallel spinor of the special type,
i.e., a parallel half-spinor, this reduces its holonomy to $\Symp(1)
\simeq\SU(2)\subset\SO(4)$.  Of course, this implies that~$(M^4,g)$ 
can be regarded as a Ricci-flat K\"ahler metric (in a $2$-parameter family
of ways, in fact).  These metrics are locally in one-to-one
correspondence with solutions of the complex Monge-Ampere equation in
two complex variables.  This has the local generality of two functions
of three variables.  The solutions are all real-analytic.

\subsubsection{Type~$(3,1)$}\label{sssec:parspin31}
Suppose~$(M^{3,1},g)$ has a nonzero parallel spinor.  Since
there is only one nonzero~$\Spin(3,1)$-orbit in~$\bbS^{3,1}\simeq\bbC^2$,
there is only one possible algebraic type of parallel spinor.
I can now apply the moving frame analysis to the coframe bundle adapted
to a single nonzero element in~$\bbS^{3,1}$.

Since the stabilizer subgroup of a nonzero vector
in~$\bbC^2$ under the action of~$\SL(2,\bbC)$ is conjugate
to the unipotent upper triangular matrices, take the
structure equations for coframes~$\omega_{i\bar\jmath}
=\overline{\omega_{j\bar\imath}}$ so that
\begin{equation}
g = \omega_{1\bar1}{\circ}\omega_{2\bar2} 
      - \omega_{2\bar1}{\circ}\omega_{1\bar2}
\end{equation} 
to have the form
\begin{equation}
d\begin{pmatrix}\omega_{1\bar1}&\omega_{1\bar2}\\ 
                \omega_{2\bar1}&\omega_{2\bar2} \end{pmatrix}
= -\begin{pmatrix}0&\alpha\\0&0\end{pmatrix}\w
   \begin{pmatrix}\omega_{1\bar1}&\omega_{1\bar2}\\ 
        \omega_{2\bar1}&\omega_{2\bar2}\end{pmatrix}
  +\begin{pmatrix}\omega_{1\bar1}&\omega_{1\bar2}\\ 
      \omega_{2\bar1}&\omega_{2\bar2}\end{pmatrix}
    \w \begin{pmatrix}0&0\\ \bar\alpha&0\end{pmatrix}.
\end{equation}
Since~$d\omega_{2\bar2}=0$, write~$\omega_{2\bar2}=dx_{2\bar2}$
for some $\bbR$-valued function~$x_{2\bar2}$.  Since $d\omega_{2\bar1}
=\omega_{2\bar2}\w\bar\alpha$, there exists locally a $\bbC$-valued 
coordinate~$x_{2\bar1}$ so that~$\omega_{2\bar1}
=dx_{2\bar1}-\bar{p}\,dx_{2\bar2}$.  
This forces $\alpha = dp + q\,dx_{2\bar2}$.  Reducing frames
to make~$p=0$ makes~$\alpha = q\,dx_{2\bar2}$ and
\begin{equation}
d\omega_{1\bar1}=-\alpha\w\omega_{2\bar1}+\omega_{1\bar2}\w\bar\alpha
    = (\bar{q}\,dx_{1\bar2}+q\,dx_{2\bar1})\w dx_{2\bar2}\,,
\end{equation}
so that there must be an $\bbR$-valued function~$f$ on an open set 
in~$\bbC\times\bbR$ so that
\begin{equation}
(\bar{q}\,dx_{1\bar2}+q\,dx_{2\bar1})\w dx_{2\bar2}
= d\bigl(f(x_{1\bar2},x_{2\bar2})\, dx_{2\bar2}\bigr).
\end{equation}
Thus, there is an $\bbR$-valued coordinate~$x_{1\bar1}$ so that
$\omega_{1\bar1} = dx_{1\bar1} 
+ f(x_{1\bar2},x_{2\bar2})\, dx_{2\bar2}$.  In particular,
the metric~$g$ is locally of the form
\begin{equation}
g = dx_{1\bar1}{\circ}dx_{2\bar2} 
      - dx_{2\bar1}{\circ}dx_{1\bar2} 
   + f(x_{1\bar2},x_{2\bar2})\, (dx_{2\bar2})^2.
\end{equation}

Conversely, via this formula, any function of $3$ variables will 
produce a $(3,1)$-metric with a parallel spinor field.  Note that $g$
will be flat if and only if the ($\bbC$-valued) curvature $2$-form
\begin{equation}
F = d\alpha 
= d\left(\frac{\p  f}{\p  x_{2\bar1}}\ dx_{2\bar2}\right)
\end{equation}
vanishes, i.e., $f$ is linear in~$x_{2\bar1}$ and~$x_{1\bar2}$.
Moreover, $g$ is Ricci-flat if and only if~$f$ 
is harmonic in the complex variable~$x_{2\bar1}$, 
which does not imply flatness. 

The conclusion is that the local Ricci-flat examples with a 
parallel spinor field depend on two (real) functions of two (real)
variables.  (The coordinate ambiguity is functions of one variable.)
Of course, this normal form is well-known in general relativity.

\subsubsection{Type~$(2,2)$}\label{sssec:parspin22}
The most interesting $4$-dimensional case, from my point of view, is that
of~$(M^{2,2},g)$ and the different possibilities for a parallel spinor.
Recall from \ref{par:spin22} that $\Spin(2,2)$ has one open 
orbit in~$\bbS^{2,2}$ and two degenerate orbits, 
which form a single~$\Pin(2,2)$ orbit.  Thus, there are two subcases:

\paragraph{Generic type}\label{par:parspin22gen}
The case of a parallel spinor field in the open orbit is very much like
that just treated.  Take the model spinor to be
\begin{equation}
s = (s_+,s_-) = \left( \begin{pmatrix} 1\\0\end{pmatrix},
                        \begin{pmatrix} 0\\1\end{pmatrix} \right).
\end{equation}
Then the tautological form~$\omega$ takes values in~$\bbR^{2,2}=\bbR(2)$
and satisfies~$d\omega = -\alpha\w\omega-\omega\w\beta$ where~$\alpha$
and~$\beta$ take values in the Lie algebra of the stabilizer of~$s_\pm$,
i.e.,
\begin{equation}
\alpha = \begin{pmatrix}0 & \alpha^1_2\\ 0 & 0\end{pmatrix}
\qquad\text{and}\qquad
\beta = \begin{pmatrix}0 & 0\\ \beta^2_1 & 0\end{pmatrix}.
\end{equation}
The structure equations then become
\begin{equation}
d\begin{pmatrix}\omega^1_1 & \omega^1_2\\ 
                 \omega^2_1 & \omega^2_2\end{pmatrix}
= -\begin{pmatrix}0 & \alpha^1_2\\ 0 & 0\end{pmatrix}\w
\begin{pmatrix}\omega^1_1 & \omega^1_2\\ 
                 \omega^2_1 & \omega^2_2\end{pmatrix}
- \begin{pmatrix}\omega^1_1 & \omega^1_2\\ 
                 \omega^2_1 & \omega^2_2\end{pmatrix}\w
      \begin{pmatrix}0 & 0\\ \beta^2_1 & 0\end{pmatrix}.
\end{equation}
Thus~$d\omega^2_2=0$, so there exists a function~$x^2_2$, unique up to 
an additive constant, so that~$\omega^2_2 = dx^2_2$.  
The equation~$d\omega^2_1 = \beta^2_1\w\omega^2_2$
then implies that there exist functions~$x^2_1$ and~$b$ on the frame bundle,
with~$x^2_1$ unique up to the addition of a function of~$x^2_2$, 
so that~$\omega^2_1 = dx^2_1 + b\,dx^2_2$.  Similarly, there exist functions
$x^1_2$ and~$a$ on the frame bundle,
with~$x^1_2$ unique up to the addition of a function of~$x^2_2$, 
so that~$\omega^1_2 = dx^1_2 - a\,dx^2_2$. Reducing frames so 
that $a = b = 0$ yields~$\omega^1_2 = dx^1_2$ and~$\omega^2_1 = dx^2_1$
and the structure equations now imply that~$\beta^2_1\w dx^2_2 = \alpha^1_2
\w dx^2_2 = 0$, so that there must exist functions~$p$ and~$q$ so that
$\alpha^1_2 = p\,dx^2_2$ and $\beta^2_1 = -q\,dx^2_2$.  The structure
equation
\begin{equation}
d\omega^1_1 = -\alpha^1_2\w \omega^2_1 +\beta^2_1\w\omega^1_2
= (p\,dx^2_1 + q\,dx^1_2)\w dx^2_2
\end{equation}
now implies that there must exist functions~$x^1_1$ and~$f$, with
$x^1_1$ unique up to the addition of a function of~$x^2_2$ so
that~$\omega^1_1 = dx^1_1 + f\,dx^2_2$.  Going back to the~$d\omega^1_1$
structure equation, this implies that the function~$f$ satisfies
\begin{equation}
\frac{\p  f}{\p  x^1_1} = 0,\quad
\frac{\p  f}{\p  x^2_1} = p,\quad\text{and}\quad
\frac{\p  f}{\p  x^1_2} = q.
\end{equation}
This analysis shows that there exist local
coordinates~$x^1_1,x^2_2,x^1_2,x^2_1$ and a function~$f$ on
an open set in~$\bbR^3$ so that
\begin{equation}
g = dx^1_1\,dx^2_2 
      - dx^2_1\,dx^1_2 
   + f(x^1_2,x^2_1,x^2_2)\, (dx^2_2)^2.
\end{equation}
Moreover, these coordinates are canonical up to functions of one variable.
This metric is flat if and only if the curvature forms
\begin{equation}
d\alpha^1_2 = d\left(\frac{\p  f}{\p  x^2_1}\right)\w dx^2_2
\qquad\text{and}\qquad
d\beta_1^2 = -d\left(\frac{\p  f}{\p  x^1_2}\right)\w dx^2_2
\end{equation}
both vanish, which can only happen if~$f$ is linear in~$x^2_1$ and~$x^1_2$.

This metric is Ricci-flat if and only if $f$ satisfies
\begin{equation}
\frac{\p ^2 f}{\p  x^2_1\p  x^1_2} = 0,
\end{equation}
so the Ricci-flat metrics with a generic parallel spinor 
depend on two functions of two variables.

\paragraph{Degenerate type}\label{par:parspin22deg}
Finally, consider the degenerate case, i.e., where the metric
has a parallel spinor field whose corresponding $\Spin(2,2)$-orbit 
is $3$-dimensional.  
Then, on the adapted frame bundle, 
the tautological form~$\omega$ takes values in~$\bbR^{2,2}=\bbR(2)$
and satisfies~$d\omega = -\alpha\w\omega-\omega\w\beta$ where~$\alpha$
and~$\beta$ take values in the Lie algebra of the stabilizer of~$s_+$,
i.e.,
\begin{equation}
\alpha = \begin{pmatrix}0 & \alpha^1_2\\ 0 & 0\end{pmatrix}
\qquad\text{and}\qquad
\beta = \begin{pmatrix}\beta^1_1 & \beta^1_2\\ 
         \beta^2_1 & -\beta^1_1\end{pmatrix}.
\end{equation}
The structure equations then become
\begin{equation}
d\begin{pmatrix}\omega^1_1 & \omega^1_2\\ 
                 \omega^2_1 & \omega^2_2\end{pmatrix}
= -\begin{pmatrix}0 & \alpha^1_2\\ 0 & 0\end{pmatrix}\w
\begin{pmatrix}\omega^1_1 & \omega^1_2\\ 
                 \omega^2_1 & \omega^2_2\end{pmatrix}
- \begin{pmatrix}\omega^1_1 & \omega^1_2\\ 
                 \omega^2_1 & \omega^2_2\end{pmatrix}\w
     \begin{pmatrix}\beta^1_1 & \beta^1_2\\ 
         \beta^2_1 & -\beta^1_1\end{pmatrix}.
\end{equation}
This implies that the form~$\omega^2_1\w\omega^2_2$ is parallel,
and, in particular, closed.  Thus, each point of~$M$ has an
open neighborhood~$U$ on which there exist functions
$x = (x_1,x_2)$ so that~$\omega^2_1\w\omega^2_2
= dx_1\w dx_2$.  One can then do a bundle reduction over~$U$ so that
$\omega^2_1 = dx_1$ and ~$\omega^2_2 = dx_2$.  The structure equations
for~$d\omega^2_i$ then imply that
\begin{equation}
0 = dx_1\w\beta^1_1 + dx_2\w\beta^2_1 = dx_1\w\beta^1_2 - dx_2\w\beta^1_1\,.
\end{equation}
By Cartan's Lemma, it follows that there exist functions~$q_1,\ldots,q_4$
so that
\begin{equation}
\begin{pmatrix}
-\beta^1_2\\ \beta^1_1\\ \beta^2_1\end{pmatrix}
= \begin{pmatrix} q_1 & q_2 \\ q_2 & q_3 \\ q_3 & q_4\end{pmatrix}
\begin{pmatrix} dx_1\\ dx_2\end{pmatrix}.
\end{equation}
Using this, it follows from the structure equations that 
\begin{equation}
d\omega^1_1 \equiv d\omega^1_2 \equiv 0 \mod dx_1,dx_2\,.
\end{equation}
Consequently, each point of~$U$ has an open neighborhood~$V\subset U$
on which there exist functions~$y = (y_1,y_2)$ for 
which~$\omega^1_i \equiv dy_i \mod dx_1,dx_2$.  Obviously,
the functions~$(x_1,x_2,y_1,y_2)$ are independent on~$V$, so
by shrinking~$V$ if necessary, one can assume that they form a cubic 
coordinate system on~$V$. The congruences above show that
\begin{equation}
g = \omega^1_1\,\omega^2_2 - \omega^1_2\,\omega^2_1
  = dy_1\,dx_2 - dy_2\,dx_1 + s^{ij}(x,y)\,dx_i\,dx_j
\end{equation}
for some functions~$s^{ij} = s^{ji}$ on the range of the coordinate
chart~$(x,y):V\to\bbR^4$.  
By a final reduction of the bundle structure over~$V$, 
one can arrange
\begin{equation}
\omega^1_1 = dy_1 + s^{12}\,dx_1 + s^{22}\,dx_2\,,
\qquad\qquad
\omega^1_2 = dy_2 - s^{11}\,dx_1 - s^{12}\,dx_2\,.
\end{equation}
On this bundle,~$\alpha^1_2 = p^1\,dx_1 + p^2\,dx_2 + r^1\,dy_1 + r^2\,dy_2$
for some functions~$p^1,p^2,r^1,r^2$.

Now, going back to the structure equations, one finds that they force
\begin{equation}
\frac{\p  s^{11}}{\p  y_2} 
= \frac{\p  s^{12}}{\p  y_1} 
\qquad\text{and}\qquad
\frac{\p  s^{12}}{\p  y_2} 
= \frac{\p  s^{22}}{\p  y_1} \,,
\end{equation}
implying that there must be a function~$f$ on the hypercube~$(x,y)(V)\subset
\bbR^4$ so that
\begin{equation}
s^{ij} = \frac{\p ^2 f}{\p  y_i\p  y_j} \,.
\end{equation}
Conversely, given any smooth function~$f$ on a domain~$D\subset\bbR^4$,
one can define~$s^{ij}$ by the above formulae and then 
the structure equations above can be solved uniquely for the 
quantities~$p$, $q$ and $r$ (it turns out that~$r\equiv 0$ anyway). 
 Consequently, the metric
\begin{equation}
g = dy_2\,dx_1 - dy_1\,dx_2 
    +\frac{\p ^2 f}{\p  y_i\p  y_j}(x,y)
      \,dx^i\,dx^j
\end{equation}
always has a parallel spinor field of degenerate type.
Thus, these metrics depend on one arbitrary function of four variables.  
(The ambiguities in the choice of coordinates are easily seen to depend
on three functions of two variables.)  By examining the curvature of
this metric for `generic' $f$, one sees that the generic such metric
does not have more than one parallel spinor field.  In fact, the holonomy
group of the generic example is equal to the full stabilizer of a
degenerate spinor, the maximum possible.

Now, about the Einstein equations:  Using the derived formulae for
$\beta^i_j$ and~$\alpha^1_2$, one computes that
\begin{equation}
d\alpha^1_2 = S(f)\,dx_1\w dx_2 + R^{ij}(f)\,dy_i\w dx_j
\end{equation}
for certain fourth order differential operators~$S$ and~$R^{ij}=R^{ji}$ 
($1\le i,j\le 2$).  The Ricci tensor of~$g$ turns out (apart 
from an overall constant factor) to be
\begin{equation}
\Ric(g) = R^{ij}(f)\,dx_i\,dx_j\,.
\end{equation} 
Thus, the metric is Ricci-flat if and only if~$f$ satisfies a system
of three fourth order quasilinear PDE.  Although I will not
give details here (anyway, a more interesting example of this sort of
calculation will be presented later during the $7$-dimensional discussion), 
this system turns out to be involutive, with the general solution
depending on two arbitrary functions of three variables, \emph{the
same generality as in the positive definite case}.  Moreover, the
generic Ricci-flat $(2,2)$-metric with a degenerate parallel spinor field
has holonomy equal to the full stabilizer of a degenerate spinor,
again, the maximum possible.

\subsection{Dimension~$5$}\label{ssec:parspin5}

Now we move into slightly less familiar territory.

\subsubsection{Type~$(5,0)$}\label{sssec:parspin50}
In the Riemannian case, $\Spin(5)=\Symp(2)$ acts transitively on the
unit sphere in~$\bbS^{5,0}\simeq\bbH^2$, so there is only one kind
of spinor, having stabilizer subgroup~$\Symp(1)$.  This~$\Symp(1)$ maps
into~$\SO(5)$ faithfully and so lies in a copy of~$\SO(4)\subset\SO(5)$.
Thus, a Riemannian $5$-manifold
with a parallel spinor is locally the product of the metric on a line and
a Ricci-flat K\"ahler metric, which reduces our problem to the 
$4$-dimensional case.

\subsubsection{Type~$(4,1)$}\label{sssec:parspin41}
This case is considerably more interesting.
Now, $\Spin(4,1)=\Symp(1,1)$ acts transitively on the
level sets of the spinor `norm'~$\nu(s) = s^*Qs$ 
in~$\bbS^{4,1}\simeq\bbH^{1,1}$ minus ~$0\in\bbH^{1,1}$.
Thus, as explained earlier, there are two essentially different
kinds of orbits:  The first corresponding to the nonzero level
sets of~$\nu$, and the second corresponding to the zero level 
sets of~$\nu$.

\paragraph{Generic type}\label{par:parspin41gen}
If the parallel spinor field has nonzero spinor norm, then it corresponds
to a spinor in~$\bbS^{4,1}$ whose stabilizer subgroup is~$\Symp(1)$.  
Looking at the spinor squaring map, this~$\Symp(1)$ maps 
into~$\SO(4,1)$ faithfully and so lies in a copy of~$\SO(4)\subset\SO(4,1)$.
Thus a metric~$g$ of this type is locally of the form
$g = -dt^2 + \bar g$, where $\bar g$ is a Ricci-flat K\"ahler metric
on a $4$-manifold, which again reduces our problem to the 
$4$-dimensional case.

\paragraph{Degenerate type}\label{par:parspin41deg}
If the parallel spinor field has vanishing spinor norm, then it corresponds
to a spinor in~$\bbS^{4,1}$ whose stabilizer subgroup is~$G_0\simeq\bbR^3$. 
I can now apply the moving frame analysis to the coframe bundle adapted
to a such a spinor, which can be assumed to be~$s_0$, 
as defined in~\S\ref{par:spin41}.

Since the stabilizer subgroup of~$s_0$ is~$G_0$, take the
structure equations for coframes~
\begin{equation}
\omega = \begin{pmatrix}
           \omega_{1} &\omega_{2}\\
           \overline{\omega_{2}} &\omega_{1}
         \end{pmatrix}
\quad\text{where $\omega_1 = \overline{\omega_1}\,$, and}\quad
\alpha =  \begin{pmatrix}
           \phi &-\phi \\
           \phi &-\phi
         \end{pmatrix}
\quad\text{where $\phi = -\overline{\phi}$\,,}
\end{equation}
with~$d\omega = -\alpha\w\omega + \omega\w\alpha^*$.  It simplifies
the calculations to set~$\omega_1 = \rho+\xi$ and~$\omega_2 = \rho+\sigma$
where~$\rho$ and~$\xi$ are $\bbR$-valued while $\sigma$ is~$\Im\bbH$-valued.
Then the structure equations are expressed as
\begin{equation}
d\xi = 0, \qquad d\sigma = -2\xi\w\phi,
\qquad d\rho = -\phi\w\sigma+\sigma\w\phi\,.
\end{equation}
Now, by the first equation, there must exist a local coordinate~$x$,
unique up to an additive constant,
so that~$\xi = dx$.  By the second equation~$d\sigma = 2\phi\w dx$,
so, locally, there exist functions~$s$ and~$h$ with 
values in~$\Im\bbH$ so that~$\sigma = ds + 2h\,dx$.  The function~$s$
is unique up to the addition of an $\Im\bbH$-valued function of~$x$.
The second equation now implies that~$\phi = dh + p\,dx$ for some
unique $\Im\bbH$-valued function~$p$.  Now reduce frames to make~$h = 0$
(which can clearly be done).  Then the structure equations so far
say that~$\xi = dx$, $\sigma = ds$, and ~$\phi = p\,dx$.  The third
structure equation now reads
\begin{equation}
d\rho = -\phi\w\sigma+\sigma\w\phi = (p\,ds + ds\,p)\w dx\,
\end{equation}
from which it follows that there exist $\bbR$-valued functions~$r$ and $f$ 
so that~$\rho = dr + f\,dx$, where~$r$ is unique up to the addition of
a function of~$x$.  The third structure equation then further implies
that
\begin{equation}
df \equiv p\,ds + ds\,p \mod dx,
\end{equation}
so that $f$ is a function of~$x$ and~$s$ (and, moreover, that $p$ is
essentially one-half the gradient of~$f$ in the~$s$ variables). 

Thus, the calculations so far have shown that any metric of
type~$(4,1)$ with a null parallel spinor field has local coordinate 
charts $(x,s,r):U\to\bbR\times\Im\bbH\times\bbR$ 
in which the metric can be written in the form
\begin{equation}
g = d\bar s\,ds - 2 dr\,dx - (1 + 2f(x,s))\,dx^2
\end{equation}
where~$f$ is an arbitrary function of four variables.  Conversely,
for any sufficiently differentiable function~$f$ of four variables,
the above formula defines a metric that has a parallel null
spinor field, since, setting~$\xi = dx$, $\sigma = ds$, 
$\rho = dr + f\,dx$, the structure equations above will be satisfied
by taking~$\phi = p\,dx$ where ~$p$ is the unique solution of the
equation~$df \equiv p\,ds + ds\,p \mod dx$.

Since coordinate charts of the above form are determined by the metric 
up to a choice of functions of one variable, the type $(4,1)$
metrics possessing a parallel null spinor field depend on 
one arbitrary function of four variables.

The metric~$g$ will be flat if and only if the 
connection form~$\phi = p\,dx$ is closed, which is the same thing 
as saying that~$f$ is linear in~$s$.  Computation shows that the Ricci 
curvature of~$g$ vanishes if and only if~$f$ 
is harmonic in the~$s$-variables.  Consequently, the Ricci-flat metrics 
of this type depend on two functions of three variables up to 
diffeomorphism, \emph{exactly as in the positive definite case}.  

\subsubsection{Type~$(3,2)$}\label{sssec:parspin32}
Since~$\Spin(3,2)\simeq\Symp(2,\bbR)$ with~$\bbS^{3,2}\simeq\bbR^4$,
the standard representation of~$\Symp(2,\bbR)$, it follows that all
of the nonzero elements of~$\bbS^{3,2}$ belong to a single $\Spin(3,2)$-orbit.
Thus, there is only one type of parallel spinor for~$(3,2)$-metrics. 
Since this is a `split' case, this orbit must be the pure spinor orbit.
Consequently, this case is treated in~\S\ref{sssec:parspin_pure_odd}, so
I will not consider it further here.

\subsection{Dimension~$6$}\label{ssec:parspin6}

In this section,
I will describe the less well-known classification of metrics with
parallel spinors in dimension~$6$ and types~$(6,0)$, $(5,1)$, and $(3,3)$.

\subsubsection{Type~$(6,0)$}\label{sssec:parspin60}
In the Riemannian case, $\Spin(6)=\SU(4)$ acts transitively on the
unit sphere in~$\bbS^{6,0}\simeq\bbC^4$, so there is only one kind
of spinor, having stabilizer subgroup~$\SU(3)$.  This~$\SU(3)$ maps
into~$\SO(6)$ as the standard representation, so a Riemannian $6$-manifold
with a parallel spinor is a Ricci-flat K\"ahler manifold.  As is well-known,
these are determined by a convex solution of the complex Monge-Ampere
equation and so depend on two functions of five variables.

\subsubsection{Type~$(5,1)$}\label{sssec:parspin51}
In the Lorentzian case, $\Spin(5,1)=\SL(2,\bbH)$ acting 
on~$\bbS^{5,1}\simeq\bbH^2\oplus\bbH^2$, has several types of
orbits, as laid out in~\S\ref{par:spin51}.  Each of these will
be treated in turn.

\paragraph{Generic type}\label{par:parspin51gen}
Suppose that the metric has a parallel spinor field whose
associated orbit in~$\bbS^{5,1}$ has dimension~$12$.  Then the 
stabilizer of an element of this orbit is isomorphic to~$\Symp(1)$
and is hence compact.  Moreover, examining the vector representation
of~$\Spin(5,1)$ on~$\bbR^{5,1}$, one sees that this~$\Symp(1)$ gets
mapped into a copy of an~$\SU(2)\subset\SO(4)$ fixing an orthogonal
$2$-plane of type~$(1,1)$.  It follows from the generalized
de~Rham splitting theorem then that the metric is a local product
of flat~$R^{1,1}$ with a $4$-dimensional Ricci-flat K\"ahler metric.

\paragraph{Null type}\label{par:parspin51nul}
Suppose next that the metric has a parallel spinor field
whose associated orbit in~$\bbS^{5,1}$ is the $11$-dimensional
null orbit.  This case is more interesting.  The stabilizer is
now four dimensional and abelian, as was described in~\S\ref{par:spin51}.

This case is formally very much like the cases treated 
in~\S\ref{sssec:parspin21},
~\S\ref{sssec:parspin31}, and ~\S\ref{par:parspin41deg}, 
so I will not go into details, but just give the results.

One shows that a $(5,1)$-metric with a parallel spinor field
of this type always has local coordinates~$x  = 
(x_{1\bar1},x_{1\bar 2},x_{2\bar2}):U\to\bbR\times\bbH\times\bbR)$
in which the metric can be written in the form
\begin{equation}
g = - dx_{1\bar1}\,dx_{2\bar2} + |dx_{1\bar2}|^2 
      - g(x_{1\bar2},x_{2\bar2})\,{dx_{2\bar2}}^2
\end{equation}
where~$g$ is a smooth function on the open set~$(x_{1\bar2},x_{2\bar2})(U)
\subset\bbH\times\bbR$.  These coordinates are unique up to a choice
of arbitrary functions of one variable.  Thus, metrics of this type depend 
on one arbitrary function of five variables.

The Ricci tensor of such a metric vanishes if and only if~$g$ is
harmonic in the~$x_{1\bar2}$ variables.  Thus, the Ricci-flat metrics
of this kind depend locally on two arbitrary functions of four variables.

\paragraph{Degenerate type}\label{par:parspin51deg}
Finally, suppose that the metric has a parallel spinor field
whose associated orbit in~$\bbS^{5,1}$ is one of the two $7$-dimensional
degenerate orbits, i.e., the spinor field is either of positive chirality
or negative chirality.  By switching orientations, it can be assumed that
the spinor is of positive chirality, so I will do this for the rest
of the discussion. 

Suppose, then, that~$(M^{5,1},g)$ is a pseudo-Riemannian manifold with
a degenerate, positive chirality parallel spinor field.
The structure equations of the adapted coframe bundle in this case,
where~$\omega = \omega^*$ takes values in quaternion Hermitian $2$-by-$2$
matrices and~$\alpha$ takes values in the Lie algebra of the stabilizer of
the standard first basis element of~$\bbS^{5,1}_+ = \bbH^2$ are
$d\omega = -\alpha\w\omega +\omega\w\alpha^*$, where
\begin{equation}
\omega = \begin{pmatrix}\omega_{1\bar1}&\omega_{1\bar2}\\
                        \omega_{2\bar1}&\omega_{2\bar2}\end{pmatrix}
\qquad\text{and}\qquad
\alpha =  \begin{pmatrix}0&\alpha^1_2\\0&\alpha^2_2\end{pmatrix}\,.
\end{equation}
and~$\alpha^2_2$ takes values in~$\Im\bbH$.

It follows from the structure equations~$\omega_{2\bar2}$
is well-defined on the manifold and is a parallel null $1$-form.
In particular, it is closed, so that, locally,
one can introduce a $\bbR$-valued function~$x_{2\bar2}$, unique up to
an additive constant, so that~$\omega_{2\bar2} = dx_{2\bar2}$.  

The $g$-dual vector field (also null) will be denoted~$E_{1\bar1}$.
The structure equations then give
\begin{equation}
d\omega_{1\bar2} = -\alpha^1_2\w dx_{2\bar2} -\omega_{1\bar2}\w\alpha^2_2\,.
\end{equation}
This equation has an interesting interpretation.  It says that, on
each ($5$-dimensional) leaf of~$dx_{2\bar2} = 0$, 
the metric~$g$ pulls back to be the positive semidefinite
quadratic form~$|\omega_{1\bar2}|^2$ and that this
restricted quadratic form is constant along its null curves, i.e.,
the integral curves of~$E_{1\bar1}$. Thus, this quadratic form is 
well-defined on the quotient of the leaf by the (parallel)
family of null geodesics defined by~$E_{1\bar1}$. 
Moreover, the quotient metric on each leaf has holonomy~$\Symp(1)\subset
\SO(4)$, i.e., it defines a Ricci-flat K\"ahler structure on
the $4$-dimensional quotient space of each~$dx_{2\bar2}$-leaf.
Geometrically, if one considers the quotient~$\bar M$ by the
$E_{1\bar1}$ curves, it locally fibers over~$\bbR$ canonically 
(up to translation) in the
form~$x_{2\bar2}:\bar M\to\bbR$ where the $4$-dimensional fibers
are Ricci-flat K\"ahler manifolds.

Pursuing the structure equations further, the equation
\begin{equation}
d\omega_{1\bar1} = -\alpha^1_2\w\omega_{2\bar1} 
                   + \omega_{1\bar2}\w{\overline{\alpha^1_2}}
\equiv 0 \mod \omega_{2\bar1},\omega_{1\bar2},dx_{2\bar2}
\end{equation}
implies that there exists a function~$x_{1\bar1}$, 
locally defined on~$M$ so that~$\omega_{1\bar1} \equiv dx_{1\bar1}
\mod \omega_{2\bar1},\omega_{1\bar2},dx_{2\bar2}$.  This function
is unique up to the addition of a function constant along the
integral curves of~$E_{1\bar1}$, i.e., a function on~$\bar M$
(i.e., a function of five variables). 

Once~$x_{1\bar1}$ has been chosen, there is a unique reduction of
the structure bundle 
for which~$\omega_{1\bar1} = dx_{1\bar1} + f\,dx_{2\bar2}$ for some
$\bbR$-valued function~$f$.  This implies
\begin{equation}
-\alpha^1_2\w\omega_{2\bar1} 
                   + \omega_{1\bar2}\w{\overline{\alpha^1_2}} 
 \equiv d\omega_{1\bar1} \equiv 0 \mod dx_{2\bar2}.
\end{equation}
In particular, this implies that~$\alpha^1_2\equiv 0 \mod
\omega_{2\bar1},\omega_{1\bar2},dx_{2\bar2}$, so that
\begin{equation}
df\w dx_{2\bar2} = -\alpha^1_2\w\omega_{2\bar1} 
                   + \omega_{1\bar2}\w{\overline{\alpha^1_2}} 
\end{equation}
implies that~$df\equiv 0 \mod \omega_{2\bar1},\omega_{1\bar2},dx_{2\bar2}$,
i.e., that ~$f$ is constant along the~$E_{1\bar1}$ curves and so
is a function on~$\bar M$.

Conversely, starting with a $1$-parameter family of Ricci-flat
manifolds~$x:\bar M\to\bbR$, one can attempt to reconstruct a $(5,1)$
metric as follows:  Locally, choose an $\bbH$-valued $1$-form
$\eta$ on~$\bar M$ so that, on each $x$-fiber, it is a section of the
associated~$\SU(2) = \Symp(1)$ coframe bundle, i.e., so that the metric
on each $x$-fiber is given by~$|\eta|^2$ and the three parallel
self-dual $2$-forms are the components of the ~$\Im\bbH$-valued 
$2$-form~$\eta\w\bar\eta$.
This determines~$\eta$ modulo~$dx$ up to right multiplication by
a function with values in the unit quaternions, i.e.,~$\Symp(1)$.

There is then an $\Im\bbH$-valued $1$-form~$\phi$, unique modulo~$dx$, 
so that~$d\eta \equiv-\eta\w\phi\mod dx$.  In other words, there exists a
$\bbH$-valued $1$-form~$\psi$ so that
\begin{equation}
d\eta = -\psi\w dx -\eta\w\phi.
\end{equation}
Consider the effect of different choices.  Let~$\eta' = \eta + p\,dx$
where~$p$ is a function with values in~$\bbH$, and let~$\phi' = \phi + q\,dx$
where~$q$ takes values in~$\Im\bbH$.  Then 
\begin{equation}
d\eta' = - \psi'\,dx - \eta'\w\phi'
\end{equation}
where~$\psi' = dp -p\,\phi +\eta\,q +\psi + r\,dx$ for some $\bbH$-valued
function~$r$.  If this system is to satisfy the structure equations
above, then it will have to satisfy
\begin{equation}
-\psi'\w\overline{\eta'}  +\eta'\w{\overline{\psi'}} \equiv 0 \mod dx,
\end{equation}
i.e., it must be possible to choose~$p$ and~$q$ so that
\begin{equation}
\Re\left((dp -p\,\phi +\eta\,q +\psi)\w\eta\right) \equiv  0 \mod dx.
\end{equation}
Rewriting this slightly, this becomes
\begin{equation}
d\left(\Re(p\eta)\right) \equiv 
 -\Re(\psi\w\eta) - \Re(\eta\,q\,\eta) \mod dx.
\end{equation}
The term~$\Re(p\eta)$ represents a $1$-form on each 
$x$-fiber and the term~$\Re(\eta\,q\,\eta)$ 
represents an arbitrary anti-self dual $2$-form on each fiber.  In
other words, the above equation represents determining the 
$1$-form~$\Re(p\eta)$ by specifying the self-dual part of its exterior
derivative.  This is, of course, an underdetermined elliptic equation
and so can always be solved locally.  

Suppose that such a solution has been found. (Actually, it is a $1$-parameter
family of such solutions, varying with~$x$.)  Once this has been done, 
the equation~$\Re(\psi'\w\overline{\eta'})\equiv0\mod dx$ is satisfied
and, then, by choosing~$r$ appropriately, one can arrange that~
$\Re(\psi'\w\overline{\eta'})=0$ (not just modulo~$dx$.  

For notational clarity, drop the primes and assume 
that~$\Re(\psi\w\overline{\eta})=0$.  Then the metric
\begin{equation}
g = -dy\circ dx + |\eta|^2 + f\,dx^2
\end{equation}
where $f$ is an arbitrary function on~$\bar M$, will satisfy the
structure equations necessary to be a metric of the desired type.
A count of the ambiguity in the construction shows that the solutions
depend on two arbitrary functions of five variables. (One is~$f$
and the other is the arbitrariness in the choice of~$p$.) 

Thus, the conclusion is that these metrics depend locally on two 
arbitrary functions of five variables.

I have not completed the analysis of the Einstein equations in 
this case, but hope to return to it in the future.

\subsubsection{Type~$(4,2)$}\label{sssec:parspin42}
In this case, as explained in~\S\ref{par:spin42},
there are two kinds of orbits.

\paragraph{Generic type}\label{par:parspin42gen}
The generic orbits in~$\bbS^{4,2}\simeq\bbC^{2,2}$ are the ones
on which the spinor norm is nonzero.  Each of these orbits
is a hypersurface and the stabilizer of a point in such a hypersurface
is a subgroup of~$\SU(2,2)$ that is conjugate to~$\SU(2,1)$.  Moreover,
in the spinor double cover, this subgroup is represented faithfully 
as a subgroup of~$\SO(4,2)$ that is conjugate to the standard~$\SU(2,1)$.
Consequently, these metrics are simply the Ricci-flat pseudo-K\"ahler
metrics of type~$(2,1)$.  In this respect, their analysis is essentially
the same as the analysis in the positive definite case.  The local
metrics of this kind depend on two functions of five variables.

\paragraph{Null type}\label{par:parspin42nul}
However, the situation changes when the spinor field is null.  Now
the subgroup of~$\SO(4,2)$ is not semi-simple, even though it is 
also of dimension~$8$.  I have not completed the analysis of this
case, so I will leave it for later.

\subsubsection{Type~$(3,3)$}\label{sssec:parspin33}
Now consider the split case, where $\Spin(3,3)\simeq\SL(4,\bbR)$ acts 
anti-diagonally on the sum of the two half-spinor subspaces 
$\bbS^{3,3}_\pm$.  

\paragraph{Generic type}\label{par:parspin33gen}
For the generic spinor orbit, the stabilizer subgroup is a copy 
of~$\SL(3,\bbR)\subset\SL(4,\bbR)$ and its action on~$\bbR^{3,3}$
is reducible, as~$\bbR^{3,3} = \bbR^3\oplus\bbR^3$, where the two
subspaces are null.  In fact, the action of~$\SL(3,\bbR)$ and
the quadratic form are just
\begin{align}
a\cdot (v_+,v_-) &= (a\,v_+,\ (a^*)^{-1}v_-)\\ 
Q(v_+,v_-) &= v^*_-\,v_+\,, 
\end{align}
for~$a\in\SL(3,\bbR)$ and~$v_\pm\in\bbR^3$.  

Consequently, it is not difficult
to show that a metric with this holonomy must have local
coordinates~$(x^i,y_j)$ in which it can be expressed in the form
\begin{equation}
g = \frac{\p ^2f}{\p  x^i\p  y_j}\,dx^i{\circ}dy_j
\end{equation}
where~$f$ satisfies the real Monge-Ampere equation
\begin{equation}
\det\left({\frac{\p ^2f}{\p  x^i\p  y_j}}\right)=1.
\end{equation}
Thus, the $(3,3)$-metrics with a generic parallel spinor depend on 
two functions of five variables, just as in the $(6,0)$ case.  
Moreover, these metrics are all Ricci-flat, just as in the $(6,0)$ case.

\paragraph{Null type}\label{par:parspin33null}
On the other hand, if the spinor is on the null orbit, the situation is
rather different.  Now the stabilizer subgroup of~$\Spin(3,3)$ is a
conjugate of the subgroup~$G$ consisting of matrices of the form 
\begin{equation}
\begin{pmatrix} 1&*&*&*\\ 0&*&*&*\\ 0&*&*&*\\ 0&0&0&1\end{pmatrix}
\end{equation}
The character of these solutions will be somewhat different. 
In the interest of time, let me just state the result, whose proof
is quite similar to the previous proofs.  One shows that a $(3,3)$
metric with a null parallel spinor of this kind always has
local coordinates~$(x^1,x^2,x^3,y_1,y_2,y_3)$ in which the metric
has the form
\begin{equation}
g = dy^i\,dx_i + f_{11}(x,y)\,(dx^1)^2 + 2\,f_{12}(x,y)\,dx^1\,dx^2
              +  f_{22}(x,y)\,(dx^2)^2
\end{equation} 
where the functions~$f_{11}$, $f_{12}=f_{21}$, and $f_{22}$ satisfy
the two \emph{constraint equations}
\begin{equation}
\frac{\p f_{11}}{\p y_1} + \frac{\p f_{21}}{\p y_2} =
\frac{\p f_{12}}{\p y_1} + \frac{\p f_{22}}{\p y_2} = 0,
\end{equation} 
and that, conversely, every metric of this form has a parallel spinor
field of this kind.  Moreover, these coordinates are unique up to 
choices that depend on five arbitrary functions of two variables.  It
follows that metrics satisfying these conditions essentially depend on
one arbitrary function of six variables.

The calculation of the Ricci tensor follows from the calculations to
be done below in~\S\ref{ssec:parspin_pure}, so I will not redo them
here.  Instead, I will simply report that the general metric of this
kind is not Ricci-flat, but that, when one imposes the Ricci-flat
condition as a system of equations, the resulting system is in involution
and the general solution depends on two arbitrary functions of five
variables, \emph{exactly as for the case of a non-null parallel spinor 
field}.

\paragraph{Degenerate type}\label{par:parspin33degen}
Finally, consider the case where the parallel spinor field is
associated to one of the most degenerate orbits, either~$\bbS^{3,3}_+$ 
or~$\bbS^{3,3}_-$ (minus the origin, of course). Now, this is
the split case and each of these orbits constitute the pure spinors.
Thus, this is a special case of the treatment 
in~\S\ref{sssec:parspin_pure_even}, so I will not consider it further
here, except to mention that, as in the previous two cases, the
Ricci-flat solutions depend on two arbitrary functions of five
variables.

\subsection{Parallel pure spinor fields}\label{ssec:parspin_pure}

As was pointed out in~\S\ref{ssec:spin_pure}, the most degenerate
orbits in the split cases~$\Spin(p{+}1,p$ and~$\Spin(p,p)$ are the 
so-called `pure' spinors.   The stabilizer of a pure spinor in
either case maps under the double covering to the stabilizer of a 
maximal null $p$-vector in~$\bbR^{p+1,p}$ or~$\bbR^{p,p}$, respectively.  
Thus, having a parallel pure spinor field (i.e., of the most degenerate
type) is the same as having a parallel null $p$-plane field.  From 
that point of view, the metrics with this property are easily 
analyzed.  Equivalent normal forms to the ones derived below 
have been derived independently by Ines Kath~\cite{iKpurenote}.
My main interest is in how this condition interacts with the Einstein
condition, which I explain at some length.

\subsubsection{The odd case}\label{sssec:parspin_pure_odd}
Suppose that~$(M^{p{+}1,p},g)$ is a metric with a parallel
null $p$-plane field.  Consider the bundle of coframes of the form
\begin{equation}
\omega = \begin{pmatrix}\zeta\\ \xi\\ \eta \end{pmatrix}
       = \begin{pmatrix}\zeta\\ \xi^i\\ \eta_i \end{pmatrix}
\end{equation}
(where lower case Latin indices range from~$1$ to~$p$ and the summation
convention will be in force) with the
property that $g = \zeta^2 + 2\eta^i\,\xi_i$ and the parallel null 
$p$-form is~$\xi = \xi^1\w\cdots\xi^p$.  The hypothesis that $\xi$
is parallel implies that the Levi-Civita connection $1$-form~$\alpha$
associated to~$\omega$ will have the form
\begin{equation}
\alpha = \begin{pmatrix}0 & -\tau^* & 0\\ 
                        0 & \phi & 0\\ 
                        \tau & \sigma & -\phi^*\end{pmatrix}
= \begin{pmatrix}0 & -\tau_j & 0\\ 
                        0 & \phi^i_j & 0\\ 
                        \tau_i & \sigma_{ij} & \phi^j_i\end{pmatrix} ,     
\end{equation}
where~$\tr\phi = 0$ and $\sigma + \sigma^* = 0$.

The first structure equation is~$d\omega = - \alpha\w\omega$.  In 
particular, this implies that~$d\xi = -\phi\w\xi$, so that there
exists (locally) a submersion~$x:U(\subset M)\to\bbR^p$ 
so that~$\xi = f^{-1}\,dx$ where~$f:U\to\SL(p,\bbR)$ is some smooth
mapping.  By an allowable change of coframe, it can be assumed that
$f \equiv \I_p$, so do this.  Thus,~$\xi^i = dx^i$, implying that
\begin{equation}
0 = d\xi^i = -\phi^i_j\w\xi^j = -\phi^i_j\w dx^j.
\end{equation}
By Cartan's Lemma, this implies that there exist 
functions~$f^i_{jk}=f^i_{kj}$ on~$U$ so that~$\phi^i_j = f^i_{jk}\,dx^k$.
Since~$\phi$ has trace equal to zero, it follows that~$f^i_{ij} = 0$.

Now, the first structure equation gives~$d\zeta = \tau_i\w\xi^i \equiv 0
\mod dx^1,\ldots,dx^p$.  Consequently, there exists a function~$z$ on~$U$
(shrunken, if necessary) so that~$\zeta = dz + t_i\,dx^i$.  By an allowable
change of coframe, it can be assumed that the~$t_i$ are all zero, so do
this.  This now implies that~$\zeta = dz$, so
\begin{equation}
0 = d\zeta = \tau_i\w dx^i,
\end{equation}
implying, again, by Cartan's Lemma, that there exist functions~$t_{ij}=t_{ji}$
so that~$\tau_i = t_{ij}\,dx^j$.  

Now, the structure equations imply that~
\begin{equation}
d\eta = -\tau\w\zeta - \sigma\w dx + \phi^*\w\eta
      \equiv 0 \mod dx^1,\ldots,dx^p
\end{equation}
so it follows that, after shrinking~$U$ if necessary, there is a 
function~$y:U\to \bbR^n$ so that $\eta\equiv dy \mod dx^1,\ldots,dx^p$.
I.e., there exist functions~$f_{ij}$ on~$U$ so that
\begin{equation}
\eta_i = dy_i + f_{ij}\,dx^j.
\end{equation}
Applying an allowable coframe change, I can arrange that~$f_{ij}
=f_{ji}$, so assume this from now on.  Substituting this formula
back into the structure equation for~$d\eta$ and using the skewsymmetry
of~$\sigma$ and the trace-free property of~$\phi$, it follows that
the functions~$f_{ij}$ must satisfy the $p$ first order equations
\begin{equation}
\frac{\p  f_{ij}}{\p  y_j} = 0.\label{eq:firstordertracecondition}
\end{equation}

Thus, it has been shown that a $(p{+}1,p)$-metric that possesses
a parallel pure spinor field has local coordinate 
charts~$(x,y,z):U\to\bbR^{2p+1}$ in which the metric can be expressed
as
\begin{equation}
g = dz^2 + 2\,dy_i\,dx^i +  2 f_{ij}(x,y,z)\,dx^i\,dx^j
\label{eq:pureparspinmetric_odd}
\end{equation}
where the functions~$f_{ij}=f_{ji}$ 
satisfy~\eqref{eq:firstordertracecondition}.

Conversely, I claim that a metric that can be written in this form
does possess a parallel pure spinor field.  To see this, it suffices
to take the coframing
\begin{equation}
\zeta = dz,\qquad\qquad 
\xi^i = dx^i,\qquad\qquad 
\eta_i = dy_i + f_{ij}\,dx^j
\end{equation}
and verify that setting
\begin{equation}
\begin{split}
\phi^i_j &= -\frac{\p  f_{jk}}{\p  y_i}\,dx^k,
\qquad\qquad \tau_i = \frac{\p  f_{ik}}{\p  z}\,dx^k,
\qquad\text{and}\\
\sigma_{ij} &= \left( \frac{\p  f_{ik}}{\p  x^j}
                    -\frac{\p  f_{jk}}{\p  x^i}
               + f_{il}\frac{\p  f_{jk}}{\p  y_l}
               - f_{jl}\frac{\p  f_{ik}}{\p  y_l}\right)\,dx^k,
\end{split}
\label{eq:connectionforms}
\end{equation}
satisfies the structure equations.  (Note that 
\eqref{eq:firstordertracecondition} is needed in order for~$\phi$ to
be trace-free.)

Thus, the $(p{+}1,p)$-metrics with
a parallel pure spinor field depend essentially on $\frac12p(p{+}1)-p
=\frac12p(p{-}1)$ arbitrary functions of~$2p{+}1$ variables.  (The
ambiguity in the choice of these coordinates is measured in 
functions of~$p$ variables, which is negligible.)

\paragraph{Curvature and holonomy}\label{par:parspin_pure_odd_curv}
I am now going to show that the metrics of this type do not, generally
have any more parallel spinors by showing that the holonomy group of the
generic metric of this kind is equal to the full stabilizer of a 
null $p$-vector.  This will be done by examining the curvature
of such a metric.

 The components of the curvature 
$2$-form~$\Theta = d\alpha+\alpha\w\alpha$ are
\begin{equation}
\begin{split}
\Phi^i_j &= d\phi^i_j + \phi^i_k\w\phi^k_j\,,\\
T_i &= d\tau_i - \phi_i^j\w\tau_j\,,\\
\Sigma_{ij} &= d\sigma_{ij}
             -\phi_i^k\w\sigma_{kj}+\sigma_{ik}\w\phi^k_j
          -\tau_i\w\tau_j\,.\\
\end{split}
\end{equation}
Note that the expressions~\eqref{eq:connectionforms}
for the components~$\phi^i_j$, $\tau_i$, and $\sigma_{ij}$ are all linear
combinations of the~$\xi^i$, i.e., of~$dx^i,\ldots,dx^p$. One consequence of 
this fact is that the curvature 2-forms must all lie in the 
ideal~$\cX$ generated by~$dx^i,\ldots,dx^p$.  

Now, let~$\eug\subset\euso(p{+}1,p)$ be the Lie algebra of the 
stabilizer of the null $p$-vector as described above.
By the Ambrose-Singer holonomy theorem, the Lie algebra~$\euh\subset\eug$ 
of the holonomy group at~$0\in \bbR^{2p+1}$ is spanned by the matrices 
of the form
\begin{equation}
P_\gamma^{-1}\,\Theta(v,w)\,P_\gamma
\end{equation}
where~$\gamma:[0,1]\to\bbR^{2p+1}$ is a differentiable curve 
with~$\gamma(0)=0$
and~$v$ and $w$ are tangent vectors to~$\bbR^{2p+1}$ at~$\gamma(1)$.  
In particular,
$\euh$ contains the subspace~$\eup$ that is spanned by matrices of
the form~$\Theta(v,w)$ where~$v$ and $w$ are tangent vectors to~$\bbR^{2p+1}$ 
at~$0$. Thus, to show that~$\euh=\eug$, it suffices to show that~$\eup$
generates~$\eug$ as a Lie algebra.

Now, let~$\cX\w\cX$ denote the span of the 
$2$-forms~$\{dx^i\w dx^j\ \vrule\ 1\le i,j\le p\}$.  
Using the given expressions for the
components of~$\theta$, the components of~$\Theta$ satisfy congruences
modulo~$\cX\w\cX$ of the form
\begin{equation}
\begin{split}
\Phi^i_j &\equiv \frac{\p^2 f_{jk}}{\p z\,\p y_i}\,dx^k\w dz
   + \frac{\p^2 f_{jk}}{\p y_l\p y_i}\,dx^k\w dy_l\,,\\
T_i &\equiv -\frac{\p^2 f_{ik}}{\p z^2 }\,dx^k\w dz
  - \frac{\p^2 f_{ik}}{\p z\,\p y_l}\,dx^k\w dy_l\,,\\
\Sigma_{ij} &\equiv 
-\Biggl(
\frac{\p^2f_{ik}}{\p z\,\p x^j}
-\frac{\p^2f_{jk}}{\p z\,\p x^i}
+f_{pi}\,\frac{\p^2f_{jk}}{\p z\,\p y_p}
-f_{pj}\,\frac{\p^2f_{ik}}{\p z\,\p y_p}\\
&\qquad\qquad\qquad\qquad\qquad
+\frac{\p f_{pi}}{\p z}\,\frac{\p f_{jk}}{\p y_p}
-\frac{\p f_{pj}}{\p z}\,\frac{\p f_{ik}}{\p y_p}
\Biggr)\,dx^k\w dz\\
&\quad -\Biggl(
\frac{\p^2f_{ik}}{\p y_l\,\p x^j}
-\frac{\p^2f_{jk}}{\p y_l\,\p x^i}
+f_{pi}\,\frac{\p^2f_{jk}}{\p y_l\,\p y_p}
-f_{pj}\,\frac{\p^2f_{ik}}{\p y_l\,\p y_p}\\
&\qquad\qquad\qquad\qquad\qquad
+\frac{\p f_{pi}}{\p y_l}\,\frac{\p f_{jk}}{\p y_p}
-\frac{\p f_{pj}}{\p y_l}\,\frac{\p f_{ik}}{\p y_p}
\Biggr)\,dx^k\w dy_l\\
\end{split}
\end{equation}

Consider now a particular solution of the form
\begin{equation}
f_{ij} = {\ts\frac12}\,h^{kl}_{ij}\,y_ky_l + {\ts\frac12}\,h_{ij}\,z^2
\end{equation}
where~$h_{ij}=h_{ji}$ and $h^{kl}_{ij}=h^{lk}_{ij}=h^{kl}_{jk}$ are
constants satisfying the condition~$h^{kl}_{kj}=0$.  This
choice satisfies the constraint equations~\eqref{eq:firstordertracecondition}
and, moreover, satisfies
\begin{equation}
\Phi^i_j \equiv  h^{il}_{jk}\,dx^k\w dy_l\,,
\qquad\text{and}\qquad
T_i \equiv -h_{ik}\,dx^k\w dw\,,
\end{equation}
the congruences being taken modulo~$\cX\w\cX$.  Moreover, the
$2$-forms~$\Sigma_{ij}$ vanish to order at least~$2$ at the origin~$x=y=z=0$.

It follows that, when the constants~$h_{ij}$
and~$h^{ij}_{kl}$ are taken sufficiently generically, 
the space~$\eup$ (and hence~$\euh$) contains all the matrices of 
the form
\begin{equation}
\begin{pmatrix}0 & -r^* & 0\\ 0 & q & 0\\ r & 0 & -q^*\end{pmatrix}
\end{equation}
with~$r\in\bbR^p$ and~$q\in\eusl(p,\bbR)$.  However, the space
of such matrices generates~$\eug$.  It follows that the
holonomy group is equal to the full stabilizer~$G\subset\SO(p{+}1,p)$,
as was desired.  

It follows, moreover, that there is an open, dense condition on the $2$-jet 
of the functions~$f_{ij}$ whose satisfaction will imply that the 
corresponding metric~$g$ will have holonomy equal to~$G$.  In particular, 
such a metric will have exactly one parallel spinor, 
which will moreover, be pure.

\paragraph{The Ricci tensor}\label{par:parspin_pure_odd_ricci}
Finally, I want to examine the conditions for such a metric to be
Ricci-flat.  A calculation shows that the formula for the
Ricci tensor of the metric~$g$ defined by~\eqref{eq:pureparspinmetric_odd} is
\begin{equation}
\Ric(g) = 2\left(
\frac{\p^2 f_{jl}}{\p z^2}
-\frac{\p^2 f_{jl}}{\p x^k\,\p y_k}
+f_{mk}\,\frac{\p^2 f_{jl}}{\p y_m\,\p y_k}
- \frac{\p f_{mj}}{\p y_k}\,\frac{\p f_{kl}}{\p y_m}
\right)\,dx^j\,dx^l.
\end{equation}
Thus, the generic such metric is not Ricci-flat.   

There remains the question of how many $(p{+}1,p)$-metrics there are that
both have a parallel pure spinor field and are Ricci-flat.  
By the above formula, this is, locally, the same as asking for the 
simultaneous solutions to the overdetermined system:
\begin{equation}
\begin{split}
\frac{\p f_{ij}}{\p y_j} &= 0,\\
\frac{\p^2 f_{jl}}{\p z^2}
-\frac{\p^2 f_{jl}}{\p x^k\,\p y_k}
+f_{mk}\,\frac{\p^2 f_{jl}}{\p y_m\,\p y_k}
- \frac{\p f_{mj}}{\p y_k}\,\frac{\p f_{kl}}{\p y_m}&=0.
\end{split}
\label{eq:Ricciflat}
\end{equation}
Fortunately, this system is involutive in Cartan's sense, so that 
local solutions are guaranteed to exist, at least in the real-analytic
category.  (See~\cite{bcggg} for a discussion of what this means.)

In fact, though, it is not necessary to invoke the Cartan-K\"ahler
theory in this case, as a direct proof can be given for the existence
of solutions to the Cauchy problem.  Here is how this can be done:
Consider functions~$a_{ij}=a_{ji}$ and~$b_{ij}=b_{ji}$ on~$\bbR^{2p}$ 
with coordinates~$x^i,y_j$ and suppose that these
functions satisfy the \emph{constraint equations}
\begin{equation}
\frac{\p a_{ij}}{\p y_j} = \frac{\p b_{ij}}{\p y_j} = 0.
\label{eq:constraints}
\end{equation}
Now consider the nonlinear initial value problem
\begin{equation}
\begin{split}
\frac{\p^2 f_{jl}}{\p z^2}
&=\frac{\p^2 f_{jl}}{\p x^k\,\p y_k}
  - f_{mk}\,\frac{\p^2 f_{jl}}{\p y_m\,\p y_k}
  + \frac{\p f_{mj}}{\p y_k}\,\frac{\p f_{kl}}{\p y_m}\,,\\
\noalign{\vskip 5pt}
f_{jl}(0,x,y) &= a_{jl}(x,y)\,,\\
\frac{\p f_{jl}}{\p z}(0,x,y) &= b_{jl}(x,y)\,.\\
\end{split}
\label{eq:RicciIVP}
\end{equation}

If~$a_{ij}$ and $b_{ij}$ are real-analytic, then the Cauchy-Kowalewski 
theorem implies that there is a unique real-analytic solution~$f_{jl}$ 
to this problem on an open neighborhood of~$\bbR^{2p}\times\{0\}$ 
in~$\bbR^{2p}\times\{0\}=\bbR^{2p+1}$. It must now be shown that the
resulting functions~$f_{jl}$ satisfy the constraint equations 
\begin{equation}
\frac{\p f_{ij}}{\p y_j} = 0.
\end{equation}
in order to know that they satisfy the system~\eqref{eq:Ricciflat}.

To show this, consider the real-analytic quantities
\begin{equation}
A_{l} = \frac{\p f_{jl}}{\p y_j}\,.
\end{equation}
Using the fact that~$f_{jl}$ satisfies \eqref{eq:RicciIVP},
one computes that
\begin{equation}
\begin{split}
\frac{\p^2 A_l}{\p z^2}
&= \frac{\p^3 f_{jl}}{\p z^2\,\p y_j}
 = \frac{\p\hfil}{\p y_j}\left(\frac{\p^2 f_{jl}}{\p z^2}\right)\\ 
&= \frac{\p\hfil}{\p y_j}\left(
\frac{\p^2 f_{jl}}{\p x^k\,\p y_k}
  - f_{mk}\,\frac{\p^2 f_{jl}}{\p y_m\,\p y_k}
  + \frac{\p f_{mj}}{\p y_k}\,\frac{\p f_{kl}}{\p y_m}
\right)\\ 
&= \frac{\p^2 A_{l}}{\p x^k\,\p y_k}
  - f_{mk}\,\frac{\p^2 A_{l}}{\p y_m\,\p y_k}
  + \frac{\p f_{kl}}{\p y_m}\,\frac{\p A_{m}}{\p y_k}\,.
\end{split}
\end{equation}
(Note the very fortunate circumstance that, in expanding this last step, 
the terms that appear that cannot be expressed in terms of the~$A_l$ cancel.
It is this cancellation that ensures that the constraint equations are
compatible with the Ricci equations.)
Thus, the~$A_l$ satisfy a linear second order system of~PDE in
Cauchy-Kowalewski form. Moreover,~$A_l$ satisfies the initial conditions
\begin{equation}
A_l(0,x,y) = \frac{\p a_{ij}}{\p y_j}(x,y) = 0\,,
\qquad\hbox{and}\qquad
\frac{\p A_l}{\p w}(0,x,y) = \frac{\p b_{ij}}{\p y_j}(x,y) = 0\,.
\end{equation}
Thus, by the uniqueness of real-analytic solutions to the 
initial value problem, it follows that~$A_l(z,x,y)=0$, 
as was to be shown.

In conclusion, it follows that the Ricci-flat $(p{+}1,p)$-metrics that 
possess a parallel pure spinor depend on $p(p{-}1)$ functions of $2p$ 
variables, locally.  Moreover, examining the discussion of curvature
and holonomy of solutions in~\S\ref{par:parspin_pure_odd_curv}, one
sees that it is possible to choose the initial data for the Cauchy
problem in such a way as to construct Ricci-flat solutions with
the full stabilizer group as holonomy.  Details are left to the reader.

\paragraph{The case $p=3$}\label{par:parspin_pure_p=3_ricci}
This analysis is particularly interesting in the case~$p=3$, as I
shall now explain.  The above argument shows that the Ricci-flat
$(4,3)$-metrics with a parallel pure spinor field depend locally
on six arbitrary functions of six variables.  This is the same generality
as that for $(4,3)$-metrics with a parallel spinor field that
is not null, since these are precisely the $(4,3)$-metrics 
whose holonomy lies in~$\G^*_2$, the stabilizer of any non-null
spinor in~$\bbS^{4,3}$, see~\cite{rBr87}.  It is interesting that,
even though the orbits of the null spinors and the non-null spinors 
have the same dimension, the condition to have a null parallel
spinor field is weaker than that of having a non-null parallel
spinor field.  However, adding in the Ricci-flat condition (which
is automatic for metrics with a non-null parallel spinor 
field) restores equality between the two cases, as far as local 
generality goes.

\paragraph{The case $p=4$}\label{par:parspin_pure_p=4_ricci}
The case $p=4$ is also worth mentioning for comparing the
case of a non-null parallel spinor field with that of a pure spinor
field.  Recall from the discussion in~\S\ref{par:spin_pure_odd_low}
that the generic $\Spin(5,4)$-orbit in~$\bbS^{5,4}\simeq\bbR^{16}$
is a quadratic hypersurface.  The stabilizer of a spinor on
such an orbit is isomorphic to~$\Spin(4,3)$ and this maps to
a copy of~$\Spin(4,3)\subset\SO(4,4)\subset\SO(5,4)$ and so
stabilizes a non-null vector in~$\bbR^{5,4}$.  In particular,
a metric with a non-null parallel spinor must locally
be a product of a $1$-dimensional factor with a metric on an
$8$-manifold with holonomy in~$\Spin(4,3)$.  In particular, such
metrics are Ricci-flat and depend locally on $12$ arbitrary
functions of~$7$ variables~\cite{rBr87}.  

In contrast, a $(5,4)$-metric with a parallel pure spinor
field does not necessarily factor and need not be Ricci-flat.
Moreover, even if one imposes the Ricci-flat condition, the
local generality of such metrics is still $12$ functions of $8$ 
variables.

\subsubsection{The even case}\label{sssec:parspin_pure_even}
The even case is very similar to the odd case, so I will just state
the results and leave the arguments to the reader.

First of all, one shows that a $(p,p)$-metric~$g$ that possesses
a parallel pure spinor field has local coordinate 
charts~$(x,y):U\to\bbR^{2p}$ in which the metric can be expressed
as
\begin{equation}
g = dy_i\,dx^i +   f_{ij}(x,y)\,dx^i\,dx^j
\label{eq:pureparspinmetric_even}
\end{equation}
where the functions~$f_{ij}=f_{ji}$ 
satisfy~\eqref{eq:firstordertracecondition}.

A calculation shows that the formula for the
Ricci tensor of the metric~$g$ defined by~\eqref{eq:pureparspinmetric_even} is
\begin{equation}
\Ric(g) = -2\left(
\frac{\p^2 f_{jl}}{\p x^k\,\p y_k}
-f_{mk}\,\frac{\p^2 f_{jl}}{\p y_m\,\p y_k}
+ \frac{\p f_{mj}}{\p y_k}\,\frac{\p f_{kl}}{\p y_m}
\right)\,dx^j\,dx^l.
\end{equation}
Thus, the generic such metric is not Ricci-flat.   An examination of
the curvature of this metric shows that the generic such metric
has holonomy equal to the stabilizer of a null~$p$-vector (and hence
has only one parallel spinor field).

Finally, the combination of the constraint 
equations~\eqref{eq:firstordertracecondition} and~$\Ric(g) = 0$
forms an involutive system, whose general solution depends $p(p{-}1)$
arbitrary functions of ~$2p{-}1$ variables.  Moreover, the general 
solution has holonomy equal to the stabilizer of a null~$p$-vector 
(and hence has only one parallel spinor field).

\subsection{$(10,1)$-metrics with a parallel null spinor field}
\label{ssec:parspin10,1_null}
In this final section, I will show that there are $(10,1)$-metrics
with parallel null spinor fields whose holonomy group is the maximum
possible, namely that of the group~$H\subset\SO(10,1)$ of 
dimension~$30$ that stabilizes a null spinor in~$\bbS^{10,1}$.
The notation of~\S\ref{sssec:spin101} will be continued in this
section.  By the analysis there, the image 
group~$\rho(H)\subset\SO^\uparrow(10,1)$ has Lie algebra
\begin{equation}
\rho'(\euh) = 
\left\{
\begin{pmatrix}
   0  & y &  0  &  \overline{\yb}^*\\
0 & 0  & 2y &  0   \\
0   & 0 &   0   &   0\\
 0 &  0 & 2\,\overline{\yb} & a_2\end{pmatrix}
\ \vrule\ 
\begin{matrix} y\in\bbR,\\\noalign{\vskip2pt} \yb\in\bbO,\\
        \noalign{\vskip2pt}\ a\in\euk_1 \end{matrix}\ \right\}.
\end{equation}
Thus, the problem devolves on understanding the structure equations of
a torsion-free $\rho(H)$-structure~$B\to M^{10,1}$ of the form
\begin{equation}
\begin{pmatrix}d\omega_1\\ d\omega_2\\ d\omega_3\\ d\bfomega\end{pmatrix}
= -\begin{pmatrix}
   0  & \psi &  0  &  {}^t\bfphi\\
0 & 0  & 2\,\psi &  0   \\
0   & 0 &   0   &   0\\
 0 &  0 & 2\,\bfphi & \bftheta \end{pmatrix}\w
\begin{pmatrix}d\omega_1\\ d\omega_2\\ d\omega_3\\ d\bfomega\end{pmatrix}
\end{equation}
where~$\bfomega$ and $\bfphi$ take values in~$\bbO$ and~$\bftheta$ takes
values in the subalgebra~$\euspin(7)\subset\eugl(\bbO)$ that consists of
the elements of the form~$a_2$ with~$a\in\euk_1$.  For such a 
$\rho(H)$-structure, the Lorentzian metric~$g = -4\,\omega_1\,\omega_3 
 +{\omega_2}^2+\bfomega\cdot\bfomega$ has a parallel null spinor and~$B$
represents the structure reduction afforded by this parallel structure.
Note that the null $1$-form~$\omega_3$ is parallel and well-defined on~$M$.
It (or, more properly, its metric dual vector field) is the square of
the parallel null spinor field.

Differentiating the Cartan structure equations yields the
first Bianchi identities:
\begin{equation}
0
= \begin{pmatrix}
   0  & \Psi &  0  &  {}^t\bfPhi\\ 
0 & 0  & 2\,\Psi &  0   \\  
0   & 0 &   0   &   0\\ 
 0 &  0 & 2\,\bfPhi & \bfTheta \end{pmatrix}\w
\begin{pmatrix}\omega_1\\  \omega_2\\  \omega_3\\  \bfomega\end{pmatrix}\,.
\end{equation}
where~$\Psi  = d\psi$, $\bfPhi = d\bfphi + \bftheta\w\bfphi$, and
$\bfTheta = d\bftheta + \bftheta\w\bftheta$.

By the second line of this system, $\Psi\w\omega_3=0$, while the
first line implies that~$\Psi\w\omega_2\equiv0\bmod\bfomega$, so there must
be functions~$p$~and~$\qb$, with values in~$\bbR$ and $\bbO$ respectively,
so that
\begin{equation}
\Psi = (p\,\omega_2 + \qb\cdot\bfomega)\w\omega_3\,.
\end{equation}
Substituting this into the first line of the system yields
\begin{equation}
{}^t\bigl(\bfPhi - \qb\,\omega_2\w\omega_3\bigr)\w\bfomega = 0,
\end{equation}
so it follows that
\begin{equation}
\bfPhi = \qb\,\omega_2\w\omega_3 + \bfsigma\w\bfomega\,,
\end{equation}
where~$\bfsigma = {}^t\bfsigma$ is some $1$-form with values in the
symmetric part of~$\eugl(\bbO)$, which will be denoted~$S^2(\bbO)$
from now on.  Substituting this last equation 
into the last line of the Bianchi identities, yields
\begin{equation}
2\,\bfsigma\w\bfomega\w\omega_3 + \bfTheta\w\bfomega = \zerob.
\end{equation}
In particular, this implies that~$\bfTheta\w\bfomega = \zerob\bmod\omega_3$,
so that~$\bfTheta \equiv\bR\bigl(\bfomega\w\bfomega\bigr)\bmod\omega_3$
where~$\bR$ is a function on~$B$ with values 
in~${\cK}\bigl(\euspin(7)\bigr)$, which is the irreducible~$\Spin(7)$
module of highest weight~$(0,2,0)$ and of (real) dimension~$168$.  (This
uses the usual calculation of the curvature tensor of $\Spin(7)$-manifolds.)
Thus, set
\begin{equation}
\bfTheta = \bR\bigl(\bfomega\w\bfomega\bigr) + 2\,\bfalpha\w\omega_3\,,
\end{equation}
where~$\bfalpha$ is a $1$-form with values in $\euspin(7)$ whose entries can
be assumed, without loss of generality, to be linear combinations of
$\omega_1$, $\omega_2$, and the components of~$\bfomega$. Substituting 
this last relation into the last line of the Bianchi identities now yields
\begin{equation}
2\,\bfsigma\w\bfomega\w\omega_3 + 2\,(\bfalpha\w\omega_3)\w\bfomega = \zerob,
\end{equation}
which is equivalent to the condition
\begin{equation}
\bfsigma\w\bfomega \equiv \bfalpha\w\bfomega \bmod\omega_3.
\end{equation}
In particular, this implies that~$\bfsigma-\bfalpha\equiv0\bmod\omega_3,
\bfomega$. Since~$\bfsigma$ and~$\bfalpha$ take values in~$S^2(\bbO)$ and 
$\euspin(7)$ respectively, which are disjoint subspaces of~$\eugl(\bbO)$, 
it follows that $\bfsigma\equiv\bfalpha\equiv0\bmod\omega_3,\bfomega$.  
In particular, neither~$\omega_1$ nor~$\omega_2$ appear in the expressions 
for~$\bfsigma$ and~$\bfalpha$. Recall that, by definition, $\omega_3$ does not 
appear in the expression for~$\bfalpha$, so $\bfalpha$ must be a linear 
combination of the components of~$\bfomega$ alone.
Now, from the above equation, it follows that
\begin{equation}
\bfsigma\w\bfomega = \bfalpha\w\bfomega + \sbold\,\omega_3\w\bfomega
\end{equation}
where~$\sbold$ takes values in~$S^2(\bbO)$.  Finally, the first line of
the Bianchi identities show that ${}^t\bfomega\w\bfalpha\w\bfomega = 0$, 
so it follows that~$\bfalpha = \ab(\bfomega)$ where~$\ab$ is a function
on~$B$ that takes values in a subspace of~$\Hom\bigl(\bbO,\euspin(7)\bigr)$ 
that is of dimension~$8\cdot 21 - 56 = 112$.  By the usual weights and
roots calculation, it follows that this subspace is irreducible, with 
highest weight~$(0,1,1)$.

To summarize, the Bianchi identities show that the curvature of a 
torsion-free~$\rho(H)$-structure~$B$ must have the form
\begin{equation}
\begin{split}
\Psi &=(p\,\omega_2 + \qb\cdot\bfomega)\w\omega_3\,,\\
\bfPhi &= \qb\,\omega_2\w\omega_3 
                 + \sbold\,\omega_3\w\bfomega + \ab(\bfomega)\w\bfomega \\
\bfTheta &= \bR\bigl(\bfomega\w\bfomega\bigr) + 2\,\ab(\bfomega)\w\omega_3
\end{split}
\end{equation}
where~$\bR$ takes values in~${\cK}\bigl(\euspin(7)\bigr)$, the 
irreducible $\Spin(7)$-representation of highest weight~$(0,2,0)$ (of
dimension~$168$), $\ab$ takes values in the irreducible 
$\Spin(7)$-representation of highest weight~$(0,1,1)$ (of dimension~$112$),
$\sbold$ takes values in~$S^2(\bbO)$ (the sum of a trivial representation
with an irreducible one of highest weight~$(0,0,2)$ and of dimension~$35$),
$\qb$ takes values in~$\bbO$, and $p$ takes values in~$\bbR$.  Thus, the
curvature space ${\cK}\bigl(\rho'(\euh)\bigr)$ has dimension~$325$.
By inspection, this curvature space passes Berger's first test (i.e., 
the generic element has the full~$\rho'(\euh)$ as its range).
Thus, a structure with the full holonomy is not ruled out by this method.

To go further in the
analysis, it will be useful to integrate the structure equations,
at least locally.  This will be done by a series of observations.

To begin, notice that, since~$d\omega_3=0$, there exists, locally, 
a function~$x_3$ on~$M$ so that $\omega_3=dx_3$.  This function is 
determined up to an additive constant, and can be defined on
any simply connected open subset~$U_0\subset M$.

Since~$d\omega_2 = -2\,\psi\w\omega_3 = -2\,\psi\w dx_3$, it follows that
any point of~$U_0$ has an open neighborhood~$U_1\subset U_0$ on which
there exists a function~$x_2$ for which~$\omega_2\w\omega_3 = dx_2\w dx_3$.
The function~$x_2$ is determined up to the addition of an arbitrary function 
of~$x_3$.  In consequence, there exists a function~$r$ 
on~$B_1 = \pi^{-1}(U_1)$ so that~$\omega_2 = dx_2 - 2r\,dx_3$.  
It now follows from the structure equation for~$d\omega_2$ 
that~$\psi\w\omega_3 = dr\w dx_3$.  
Consequently, there is a function~$f$ on~$B_1$ so that~$\psi = dr + f\,dx_3$.
Since~$\Psi = d\psi$ is $\pi$-basic, it follows that~$df\w dx_3$ is 
well-defined on~$U_1$.  Consequently,~$f$ is well-defined on~$U_1$
up to the addition of an arbitrary function of~$x_3$. 

Now, since
\begin{equation}
d\omega_1 = -\psi\w\omega_2 -{}^t\bfphi\w\bfomega 
          = -(dr + f\,dx_3)\w(dx_2 - 2r\,dx_3) -{}^t\bfphi\w\bfomega,
\end{equation}
it follows that
\begin{equation}
d(\omega_1 + r\,dx_2 - r^2\,dx_3) = f\,dx_2\w dx_3 -{}^t\bfphi\w\bfomega.
\end{equation}
The fact that the 2-form on the right hand side is closed, together with
the fact that the system~$I$ of dimension~9 spanned by $dx_3$ and the 
components of~$\bfomega$ is integrable (which follows from the structure
equations), implies that there are functions~$G$ and~$\bF$ on~$B$ so
that
\begin{equation}
d(\omega_1 + r\,dx_2 - r^2\,dx_3) = d( G\,dx_3 - {}^t\bF\,\bfomega ),
\end{equation}
from which it follows that there is a function~$x_1$ on~$B$ so that
\begin{equation}
\omega_1  = dx_1 - r\,dx_2 + r^2\,dx_3 + G\,dx_3 - {}^t\bF\,\bfomega\,.
\end{equation}
The function~$x_1$ is determined (once the choices of~$x_3$ and~$x_2$ are
made) up to an additive function that is constant on the leaves of the 
system~$I$, i.e., up to the addition of an (arbitrary) function of $9$
variables.  Expanding $d( G\,dx_3 - {}^t\bF\,\bfomega )
= f\,dx_2\w dx_3 -{}^t\bfphi\w\bfomega$ via the structure equations and 
reducing modulo~$dx_3$ yields
\begin{equation}
{}^t\bigl(d\bF + \bftheta\,\bF)\w\bfomega 
\equiv {}^t\bfphi\w\bfomega \bmod dx_3\,.
\end{equation}
so that there must exist functions $\bH$ and $\ub = {}^t\ub$ so that
\begin{equation}
\bfphi = d\bF + \bftheta\,\bF + \bH\,dx_3 + \ub\,\bfomega\,.
\end{equation}
Substituting this back into the relation $d( G\,dx_3 - {}^t\bF\,\bfomega )
= f\,dx_2\w dx_3 -{}^t\bfphi\w\bfomega$ yields
\begin{equation}
dG+2\,{}^t\bF\,d\bF-{}^t\bigl(\bH-2\ub\,\bF\bigr)\,\bfomega 
\equiv f\,dx_2\bmod dx_3\,.
\end{equation}
Setting $G = g - \bF\cdot\bF$ and $\hb = \bH-2\ub\,\bF$, this becomes
\begin{equation}
dg \equiv f\,dx_2 + {}^t\hb\,\bfomega \bmod dx_3\,,
\end{equation}
with the formulae
\begin{equation}
\begin{split}
\omega_1
 &=dx_1-r\,dx_2+r^2\,dx_3+(g{-}\bF\cdot\bF)\,dx_3-{}^t\bF\,\bfomega\,,\\
\bfphi &= d\bF + \bftheta\,\bF + (\hb+2\ub\,\bF)\,dx_3 + \ub\,\bfomega\,.\\
\end{split}
\end{equation}
Now the final structure equation becomes
\begin{equation}
d\bfomega = - 2\bigl(d\bF + \bftheta\,\bF + \ub\,\bfomega\bigr)\w dx_3
 - \bftheta\w\bfomega
\end{equation}
which can be rearranged to give
\begin{equation}
d\bigl(\bfomega + 2\bF\,dx_3\bigr) 
= -\bigl(\bftheta - 2\ub\,dx_3\bigr)\w \bigl(\bfomega + 2\bF\,dx_3\bigr)\,.
\end{equation}

This suggests setting $\bfeta = \bfomega + 2\bF\,dx_3$ and writing the 
formulae found so far as
\begin{equation}
\begin{split}
\omega_1&=dx_1-r\,dx_2+r^2\,dx_3+(g{+}\bF\cdot\bF)\,dx_3-{}^t\bF\,\bfeta\,,\\
\omega_2&=dx_2 - 2r\,dx_3\,,\\
\omega_3&=dx_3\,,\cr
\bfomega &= - 2\bF\,dx_3 + \bfeta\,,\\
\noalign{\vskip5pt}
\psi &= dr + f\,dx_3\,,\\
\bfphi &= d\bF + \bftheta\,\bF + \hb\,dx_3 + \ub\,\bfeta\,,\\
\noalign{\vskip5pt}
dg &\equiv f\,dx_2 + {}^t\hb\,\bfeta \bmod dx_3\,,\\
d\bfeta &= -\bigl(\bftheta - 2\ub\,dx_3\bigr)\w\bfeta \,.
\end{split}
\end{equation}
where, in these equations, $\bftheta$ takes values in~$\euspin(7)$ and
$\ub={}^t\ub$.  Note that
\begin{equation}
-4\,\omega_1\,\omega_3 + {\omega_2}^2 +\bfomega\cdot\bfomega
= -4\,dx_1\,dx_3 + {dx_2}^2 - 4g\,{dx_3}^2 + \bfeta\cdot\bfeta.
\end{equation}

I now want to describe how these formulae give a recipe for
writing down all of the solutions to our problem.  

By the last of the
structure equations, the eight components of~$\bfeta$ describe an
integrable system of rank~$8$ that is (locally) defined on the
original $11$-manifold.  Let us restrict to a neighborhood where the
leaf space of~$\bfeta$ is simple, i.e., is a smooth manifold~$K^8$.
The equation~$d\bfeta = -\bigl(\bftheta - 2\ub\,dx_3\bigr)\w\bfeta$
shows that on~$\bbR\times K^8$, with coordinate~$x_3$ on the first
factor, there is a $\{1\}\times\Spin(7)$-structure, which
can be thought of as a $1$-parameter family of torsion-free 
$\Spin(7)$-structures on~$K^8$ (the parameter is~$x_3$, of course).  

This $1$-parameter family is not arbitrary because the matrix~$\ub$ is
symmetric.  This condition is equivalent to saying that if~$\Phi$
is the canonical $\Spin(7)$-invariant $4$-form (depending on~$x_3$, of
course) then 
\begin{equation}
{{\p \Phi}\over{\p x_3}} = \lambda\,\Phi + \Upsilon
\end{equation}
for some function~$\lambda$ on~$\bbR\times K^8$ and~$\Upsilon$ is
an anti-self dual $4$-form (via the $x_3$-dependent metric on the
fibers of~$\bbR\times K\to \bbR$, of course).  It is not hard to see
that this is seven equations on the variation of torsion-free 
$\Spin(7)$-structures and that, moreover, given any $1$-parameter variation of
torsion-free $\Spin(7)$-structures, one can (locally) gauge this family by 
diffeomorphisms preserving the fibers of~$\bbR\times K\to \bbR$ so that
it satisfies these equations.  (In fact, if $K$ is compact and the
cohomology class of~$\Phi$ in~$H^4(K,\bbR)$ is independent of~$x_3$ then
this can be done globally.)  Call such a variation {\it conformally 
anti-self dual}.

Now from the above calculations, this process can be reversed:  
One starts with any conformally anti-self dual variation
of $\Spin(7)$-structures on~$K^8$, then on $\bbR^3\times K$ one forms the 
Lorentzian metric
\begin{equation}
ds^2 = -4\,dx_1\,dx_3 + {dx_2}^2 - 4g\,{dx_3}^2 + \bfeta\cdot\bfeta
\end{equation}
where~$g$ is any function on~$\bbR^3\times K$ that satisfies ${\p  g}/
{\p  x_1} = 0$ and~$\bfeta\cdot\bfeta$ is the $x_3$-dependent metric
associated to the variation of~$\Spin(7)$-structures.  Then this Lorentzian 
metric has a parallel null spinor.  For generic choice of the variation 
of~$\Spin(7)$-structures and the function~$g$, 
this will yield a Lorentzian metric whose holonomy is the desired 
stabilizer group of dimension~30.  This can be seen by combining
the standard generality result~\cite{rBr87} 
for~$\Spin(7)$-metrics on $8$-manifolds,
which shows that for generic choices as above the curvature tensor has
range equal to the full~$\rho'(\euh)$ at the generic point, 
with the Ambrose-Singer holonomy theorem, which implies that such 
a metric will have its holonomy equal to the full group of dimension~$30$.

In particular, it follows that, up to diffeomorphism, the local solutions to 
this problem depend on one arbitrary function of $10$ variables.  One
can show that such a solution is not, in general, Ricci-flat,
in contrast to the case where a $(10,1)$-metric has a non-null parallel
spinor field.

Note, by the way, that the $4$-form~$\Phi$ will not generally be closed,
let alone parallel.  However, the $5$-form~$dx_3\w\Phi$ will be closed 
and parallel.  The other non-trivial parallel forms are the $1$-form~$dx_3$, 
the $2$-form~$dx_2\w dx_3$, and the $6$-, $9$-, and $10$-forms that are the
duals of these.

\end{document}